\documentclass[final,5p,times,twocolumn]{elsarticle}

\usepackage[utf8]{inputenc}
\usepackage[T1]{fontenc}

\usepackage{amsmath}
\usepackage{amssymb}
\usepackage{gensymb}
\usepackage[dvipsnames,svgnames,x11names]{xcolor}

\newcommand*\diff{\mathop{}\!\mathrm{d}}

\usepackage{tabularx}
\usepackage{subcaption}
\usepackage{hyperref}
\usepackage[switch, modulo, displaymath, mathlines]{lineno}
\newcommand*\patchAmsMathEnvironmentForLineno[1]{
  \expandafter\let\csname old#1\expandafter\endcsname\csname #1\endcsname
  \expandafter\let\csname oldend#1\expandafter\endcsname\csname end#1\endcsname
  \renewenvironment{#1}
  {\linenomath\csname old#1\endcsname}
  {\csname oldend#1\endcsname\endlinenomath}}
  \newcommand*\patchBothAmsMathEnvironmentsForLineno[1]{
  \patchAmsMathEnvironmentForLineno{#1}
  \patchAmsMathEnvironmentForLineno{#1*}}
  \AtBeginDocument{
  \patchBothAmsMathEnvironmentsForLineno{equation}
  \patchBothAmsMathEnvironmentsForLineno{align}
  \patchBothAmsMathEnvironmentsForLineno{flalign}
  \patchBothAmsMathEnvironmentsForLineno{alignat}
  \patchBothAmsMathEnvironmentsForLineno{gather}
  \patchBothAmsMathEnvironmentsForLineno{multline}
}
\hypersetup{
	colorlinks=true,
	linkcolor=blue,
	filecolor=magenta,      
	urlcolor=cyan,
}
\usepackage[draft, deletedmarkup=xout]{changes}
\setaddedmarkup{\textcolor{blue}{#1}}
\setdeletedmarkup{\xout{\textcolor{red}{#1}}}

\usepackage{placeins} % for \FloatBarrier

\journal{Elsevier}

\begin{document}

\begin{frontmatter}



%\title{Potential field-based multi-UAV trajectory planning for 3D visual inspection}
\title{Multi-UAV trajectory planning for 3D visual inspection of complex structures}


\author[riteh]{Stefan Ivić}
\ead{stefan.ivic@riteh.hr}
\author[mathri]{Bojan Crnković}
\ead{bojan.crnkovic@uniri.hr}
\author[riteh,cnrm]{Luka Grbčić}
\ead{luka.grbcic@riteh.hr}
\author[sdu]{Lea Matleković}
\ead{matlekovic@imada.sdu.dk}

\affiliation[riteh]{organization={Faculty of Engineering, University of Rijeka},%Department and Organization
	addressline={Vukovarska 58}, 
	city={Rijeka},
	postcode={51000}, 
	country={Croatia}}
\affiliation[mathri]{organization={Department of Mathematics, University of Rijeka},
	addressline={Slavka Krautzeka 75}, 
    city={Rijeka},
    postcode={51000}, 
    country={Croatia}}
\affiliation[sdu]{organization={Department of Mathematics and Computer Science, University of Southern Denmark},
	addressline={Campusvej 55}, 
	city={Odense},
	postcode={5230}, 
	country={Denmark}}
\affiliation[cnrm]{organization={Center for Advanced Computing and Modelling, University of Rijeka},%Department and Organization
	addressline={Radmile Matejcic 52}, 
	city={Rijeka},
	postcode={51000}, 
	country={Croatia}}

\begin{abstract}
The application of autonomous UAVs to infrastructure inspection tasks provides benefits in terms of operation time reduction, safety, and cost-effectiveness.
This paper presents trajectory planning for three-dimensional autonomous multi-UAV volume coverage and visual inspection of infrastructure based on the Heat Equation Driven Area Coverage (HEDAC) algorithm. 
The method generates trajectories using a potential field and implements distance fields to prevent collisions and to determine UAVs' camera orientation. 
It successfully achieves coverage during the visual inspection of complex structures such as a wind turbine and a bridge, outperforming a state-of-the-art method by allowing more surface area to be inspected under the same conditions. 
The presented trajectory planning method offers flexibility in various setup parameters and is applicable to real-world inspection tasks. 
Conclusively, the proposed methodology could potentially be applied to different autonomous UAV tasks, or even utilized as a UAV motion control method if its computational efficiency is improved.
\end{abstract}


\end{frontmatter}


\section{Introduction}
\label{sec:sample1}  

In recent years and with the constant demand for the automatization of processes, there has been a rising interest in research of autonomous Unmanned Aerial Vehicles (UAVs) and autonomous robotic systems in general. Known for their versatile movements and ability to explore the unstructured environment, UAVs are being used for numerous civil applications such as transportation \cite{moshref2021applications, biccici2021approach}, monitoring and surveying \cite{zhao2021structural, kim2019remote}, search and rescue \cite{martinez2021search, cho2021coverage} as well as infrastructure inspection \cite{falorca2021new}. With an increase in the UAV's autonomy, UAV applications supporting construction and infrastructure inspection are expected to dominate the total UAV market with a share of around 45\% \cite{shakhatreh2019unmanned}. Manual UAV inspections have already enhanced safety as inspectors are not challenged with working at high altitudes and with heavy equipment anymore. However, inspection tasks still require skilled operating teams and can take a long time to coordinate and execute. By raising the level of autonomy, UAVs have high potential to reduce inspection time and human labor, therefore decreasing the total inspection costs. Deploying autonomous UAVs could further improve safety as with faster and cheaper inspections, the civil infrastructure could be better maintained resulting in the prevention of accidents.  

Tackling the challenge of making UAVs autonomous requires an interdisciplinary approach, collaboration, and knowledge distribution across the fields. The agents must be equipped with accurate sensor technologies that enable the control system to direct a vehicle according to the developed path or trajectory planning strategy.
Path and trajectory planning problems are crucial to solve in order to fly autonomously. Solving them imposes additional complexities as the approach usually depends on a specific application.
Path planning refers to finding a geometric, collision-free path that achieves the mission’s objectives, while trajectory planning refers to assigning a time law to the geometric path, therefore, providing control inputs for desired movements. In most cases, path planning precedes trajectory planning; however, these two phases are not necessarily distinct \cite{gasparetto2015path} as in the case where trajectories are generated from the initial to the final position and both problems are solved at the same time. Although many studies have been published concerning path and trajectory planning for autonomous vehicles, there is a lack of research attention to multi-UAV cooperation for construction and infrastructure inspection applications \cite{shakhatreh2019unmanned}. Cooperation involves the distribution of inspection tasks and collaborative, collision-free path and trajectory planning. Cooperative planning could provide wider inspection scope, higher error tolerance, and faster completion time, but increases algorithmic complexity. Current research regarding inspection path and trajectory planning focuses on solving the Coverage Path Planning (CPP) problem which aims at determining a geometric, collision-free path covering the area or volume of interest \cite{galceran2013survey}. Afterwards, paths are usually endowed with simple point-to-point trajectories. Numerous existing publications tackle two-dimensional coverage problems, but there has been a growing need for solutions to three-dimensional coverage problems.

This paper aims to fill in the research gap within multi-UAV cooperation for infrastructure inspection by presenting a novel multi-UAV trajectory planning algorithm for three-dimensional coverage as well as applying it for the visual inspection of complex structures that usually require frequent maintenance. The algorithm takes into account vehicle parameters and environmental constraints producing somewhat smooth, collision-free trajectories. It is based on a potential field approach, designing the field using a modified Heat Equation Driven Area Coverage (HEDAC) algorithm \cite{ivic2016ergodicity} which minimizes the difference between the desired and achieved coverage, therefore leading agents through less covered space. The most important aspect when visually inspecting infrastructure employing autonomous UAV is assuring that all surfaces of the structure are included in the sensor’s field of view at least once. Therefore, we upgrade the HEDAC trajectory generation with a camera direction control which always directs the camera to the closest point of the structure. As trajectories are produced in the potential field, we leave out the geometric path determination and solve both coverage path and trajectory planning problems at the same time. When accompanied by a robust trajectory-following control onboard a UAV, HEDAC can be applied to real-world infrastructure inspection. Due to the prolonged computational time, the algorithm is unable to run onboard in real time and serve as a motion controller.

This article is organized as follows. Section~\ref{sec:two} presents an overview of the state-of-the-art research concerning three-dimensional coverage path and trajectory planning as well as the HEDAC algorithm foundation and applications to two-dimensional coverage problems. In Section~\ref{sec:three}, the adaptation of the HEDAC algorithm for three-dimensional domains is presented with the implementation based on the finite element method (FEM) and applied to a 3D unit cube ergodic coverage scenario. In Section~\ref{sec:four}, we adapt the algorithm for 3D visual inspection by introducing the distance field and the camera orientation. Furthermore, we provide a simplified surface inspection evaluation and present simulation results for the inspection of a portal, wind turbine, and bridge structures. A thorough analysis of the performance, limitations, and validation of the HEDAC approach has been conducted and presented in Section~\ref{sec:five}. Finally, we conclude in Section~\ref{sec:six} by summarizing the achievements and limitations as well as proposing future improvements.

\section{State-of-the-art}
\label{sec:two}
An autonomous infrastructure inspection requires a UAV to visually capture all elements of the structure of interest. To accomplish that, vehicles need to be equipped with an appropriate camera and traverse paths that allow the camera's field of view to cover the whole structure, ensuring desired visibility. Such paths are found by solving the CPP problem. Many algorithms have been proposed for covering two-dimensional areas for inspection and monitoring purposes, making it a well-researched topic \cite{torres2016coverage}. In recent years, more and more attention is being dedicated to three-dimensional CPP. CPP approaches are categorized as either model-based or non-model-based \cite{almadhoun2019survey}. The non-model-based approaches do not require a known reference model, therefore, are applied for the unknown environment exploration or 3D model reconstruction. The model-based approaches use a reference model to generate coverage paths, making them applicable for structural inspection \cite{tan2021automatic}.

This article focuses on producing feasible, three-dimensional coverage trajectories based on known models for structural inspections. As our approach solves the CPP problem, we present state-of-the-art achievements in the model-based 3D CPP for infrastructure inspection purposes. The literature review is divided according to the number of agents employed for inspection. Furthermore, we provide an overview of the existing research regarding the HEDAC algorithm and its applications in two-dimensional domains.

\subsection{Single-agent coverage path planning}
In the attempt to solve the CPP problem for three-dimensional structure inspection, many proposed solutions generate viewpoints around a known 3D structure and employ different techniques to find coverage paths. Viewpoints are usually generated using sampling-based techniques or grid decomposition and the coverage path is found by solving the Travelling Salesman Problem (TSP) or the Vehicle Routing Problem (VRP). Multiple optimizations are often proposed for improving the overall inspection performance. A co-optimal CPP method that simultaneously optimizes the UAV path and quality of captured images and reduces the computational complexity of the solver, all while adhering to safety and inspection requirements is presented in \cite{shang2020co}. The path optimization algorithm utilizes a Particle Swarm Optimization (PSO) framework which iteratively optimizes the coverage paths based on sampled viewpoints. The core of the method consists of a cost function that measures the quality and efficiency of a coverage inspection path and the greedy heuristic for optimization enhancement by aggressively exploring the search space. Flight trajectories are generated using global B-spline curve interpolation and the method is successfully applied for the inspection of a building, a statue, and a solar plant. Another approach where a variant of PSO is used for solving TSP is presented in \cite{phung2017enhanced}. Inspection viewpoints are generated based on a grid decomposition and performance is improved by using deterministic initialization, random mutation, and edge exchange. The focus of the article is on reducing the computation time by taking advantage of parallel computing on a GPU-based framework. Flight trajectories are not proposed, however, based on the visuals given for room and bridge inspection applications, a point-to-point controller could be employed to traverse determined paths. An alternative approach proposed in \cite{cao2020hierarchical} solves TSP for sampled viewpoints by using Google OR-Tools at two different levels to provide a path for aircraft and bridge inspections. Firstly, a high-level algorithm separates the environment into multiple subspaces at different resolutions and solves the global TSP to determine the order of visiting subspaces. Secondly, a low-level TSP solution finds paths within the subspaces for detailed coverage and resolves collisions by adding new viewpoints. A comparison with greedy TSP solver is presented in \cite{shi2021inspection} where a traversal algorithm is developed. Authors uniformly sample a point cloud model and create a bounding box around it. Viewpoints are generated by projecting sampled points onto the planes of the bounding box and a traversal path search algorithm is proposed to find paths. Paths are optimized in terms of flight distance and smoothness while considering safety constraints and sensor limitations. The proposed method is evaluated for use cases of bridge and power pylon inspection where it outperformed greedy TSP by finding simpler paths in less time. A similar approach is presented in \cite{jung2018multi}. The authors divide a structure’s volumetric map into several layers, and in each layer, a set of normal vectors of each voxel’s center point is calculated. The opposing vectors are used as viewpoints and the TSP is solved using the Lin–Kernighan Heuristic (LKH) solver. Afterward, all paths in each layer are combined to form a complete path. To quantify the inspection performance, the authors propose a coverage measure as a percentage of voxels included in the camera’s field of view, achieving the result of 99.8\%. Most of the other authors fail to provide success metrics for the determined coverage paths. 

Some single-agent CPP approaches do not have a clear boundary between viewpoint sampling and path determination. Instead, these two steps are interchangeably applied to optimize both coverage and paths. A practical application of a two-step optimization approach has been demonstrated in \cite{bircher2015structural} by using a trolley mesh model. The algorithm initially samples viewpoints based on the triangles in the mesh model and solves TSP to obtain an initial tour. Afterward, viewpoints are resampled such that final paths are shorter and full coverage is provided. Full coverage is proven by reconstructing a 3D model of the structure of interest. Authors later extend their work \cite{bircher2016three} by planning the coverage path for the wind turbine, staircase, mountain, and cone models. As a part of their work on an online exploration path planner, the same authors developed an inspection CPP algorithm based on the receding horizon approach \cite{bircher2018receding}. When employed for the inspection, the algorithm samples the next best view in a random geometric tree intending to include the structure in the camera’s field of view. Feasibility is demonstrated in a bridge inspection case and practical applicability in a real-life scenario of a simple object inspection. In Section \ref{sec:five}, we compare our approach for the bridge inspection with the receding horizon planner.

The main limitation of the described existing approaches for three-dimensional coverage path planning is the number of agents deployed for the inspection. They have considered only a single-agent inspection with a focus on the computational time to compute paths. Taking into account that these model-based planners compute paths offline, decreasing the computation time would not highly affect the inspection time. Actual inspection time seems to have more impact in achieving faster and more frequent infrastructure inspections. We intend to focus on collaborative multi-agent inspection path planning as a higher number of agents is expected to decrease the inspection time. With our approach, it is still possible to plan a single-agent trajectory in case the structure is simple enough and there is no need to reduce the inspection time. 

\subsection{Multi-agent coverage path planning}

The problem of coverage path planning is enhanced by involving multiple agents which are expected to accelerate the mission accomplishment. Many known problems such as TSP or VRP are being modified and adapted for multi-agent coverage problems \cite{almadhoun2019survey}, but there is still a lack of solutions for specific applications such as 3D infrastructure inspection. Some existing solutions for multi-agent 3D area coverage such as \cite{perez2016multi}, which divides the area in hexagons and produces a lawnmower pattern to cover it, can not be adapted for structural inspection as the focus is on covering only the top part of the 3D map. These solutions do not consider complex structures which have to be covered from all sides. A similar approach of dividing a three-dimensional model in order to determine coverage paths for multiple UAVs is presented in \cite{mansouri2018cooperative}. Flight waypoints are generated based on the horizontally sliced 3D model and converted into position-velocity-yaw trajectories. Although the approach is applied for the inspection of a wind turbine and an outdoor structure, it would not be directly applicable for more complex structures as the slicing approach would not produce appropriate paths. An example of the infrastructure inspection solution which follows the pattern of most of the previously presented single-agent approaches is found in \cite{jing2020multi}. Authors use incremental sampling to create a roadmap with information on topology, coverage, and path lengths. The multi-agent CPP is formulated as a combination of two NP-hard problems: the Set Covering Problem and the VRP problem. The problem is solved by a modified Biased Random Key Genetic Algorithm which optimizes the inspection paths for multiple agents. The issue of this and similar approaches solving a variant of a TSP is that trajectories endowed through calculated waypoints are not continuous, but rather interrupted with sharp turns. Such solutions would require UAVs to fly from point to point, stopping and turning each time a path segment is accomplished, which would result in prolonged inspection time. Alternatively, if one attempts to generate trajectories through the determined viewpoints, there is a risk that no feasible solutions for a selected vehicle model would be found. The complexity of the solution drastically increases when constraints for obstacle avoidance are introduced. Our solution addresses both CPP and trajectory planning problems at the same time by producing continuous and feasible trajectories based on a potential field approach instead of solving a multi-TSP. Furthermore, the approach is highly adaptable to different, complex structures, which we demonstrate in the upcoming sections.

\subsection{HEDAC multi-agent 2D trajectory planning}
The ergodic multi-agent area coverage method HEDAC is presented for the first time in \cite{ivic2016ergodicity}. The algorithm designs a potential field based on a steady-state heat equation using a source term that depends on the difference between the given goal density and the current coverage density. The agents' movements are directed by the gradient of that potential field, leading them to the area of interest and producing paths to achieve the goal coverage density while avoiding collisions between agents using a built-in local cooling mechanism. The algorithm is demonstrated for two-dimensional test cases and compared with the Spectral Multiscale Coverage (SMC) algorithm showing superiority in achieved coverage and computational time. Simply stated, the algorithm leads agents to attractive areas that become less attractive as the area is covered, making the approach applicable for monitoring, surveillance, and inspection purposes.

The HEDAC algorithm application is further demonstrated for autonomous non-uniform multi-agent spraying tasks and presented in \cite{ivic2019autonomous}. The algorithm combines the spraying and Dubins motion models to produce the path covering the area with the desired spraying density. The method is tested in comparison to the Lawnmower and SMC algorithms for simple geometries, as well as a realistic crop spraying case, outperforming both methods in convergence time while producing spraying density of satisfying accuracy and using less spraying media when compared with conventional spraying.

The adaptation of the HEDAC algorithm for multi-agent area search in uncertain conditions is presented in \cite{ivic2020motion}. The undetected target probability-density field is computed from the initial probability field, already achieved agents' trajectories, diverse motion, and sensing parameters of each agent. The HEDAC method directs agents toward regions with a higher concentration of undetected targets, resulting in the maximization of the target detection rate. The target search using HEDAC was simulated and compared with Lawnmower, SMC, and Receding Horizon Control (RHC) approaches, demonstrating a shorter search time to achieve the same target detection rate. 

Improvement of the HEDAC algorithm for multi-agent surveying constrained motion control in irregular domains is proposed in \cite{ivic2022constrained}. The previous implementation has considered rectangular domains due to using the finite difference method for solving a partial differential equation regulating the potential. The improved solution uses the finite element method to enable a simple and elegant application of boundary conditions and modeling of static obstacles in the domain. Optimization for collision avoidance maneuvers and the ability to set path curvature constraints improved the HEDAC algorithm's suitability for real-world 2D area surveying applications. The algorithm is tested in a synthetic scenario, as well as two realistic scenarios, providing results with high computational efficiency that enables real-time execution.

Due to the excellent coverage properties, the HEDAC approach for two-dimensional domains has been thoroughly researched and applied to different problems. The multi-agent coverage control presented in \cite{zheng2022distributed} utilizes HEDAC independently on the multiple regions of a domain partitioned using Voronoi tessellation. By achieving successful coverage, this work indicates the possibility of decentralized use of HEDAC control. An experimental utilization of HEDAC control, using a 7-axis Franka Emika robot, is employed in \cite{low2022drozbot} for the ergodic-based portrait drawing. Unlike in other HEDAC implementations, a non-stationary heat equation is used to produce the trajectories tracked by the robot-driven ink pen.  

Within this article, we present the extension of the HEDAC and its application to a three-dimensional coverage problem, specifically infrastructure inspection.

\section{Three-dimensional ergodic coverage}
\label{sec:three}
A multidimensional ergodic coverage with the HEDAC method has already been theoretically presented in  \cite{ivic2016ergodicity}, but all published improvements have been implemented and applied only to two-dimensional problems. In this section, we describe the trajectory planning approach using HEDAC and its implementation for a swarm of $n$ UAVs operating in a three-dimensional domain $\Omega \in \mathbb{R}^3$ encompassed with the boundary $\Gamma$. The algorithm designs a potential field $\psi$ that attracts agents into the uncovered space until the target coverage density $\mu_0$ is achieved. The resulting motion is described with collision-free trajectories $\mathbf{y}_i : t \to \Omega$, where $t$ is time and $i=1, \ldots, n$ are indices of UAVs. Finally, we demonstrate the approach in the crowded domain using 100 agents to successfully cover half of a unit cube while avoiding collisions. 

\subsection{Space coverage}

Three-dimensional coverage can be considered as a continuous action of a UAV along its trajectory. Mathematically, this can be defined as a convolution of instantaneous action $\phi$ and the trajectory, which results in a field occupying the space around the path of a UAV. For practical reasons, we use a radial basis function, precisely the three-dimensional Gaussian function, as the instantaneous action. The smoothness of the Gaussian function allows better stability and smooth motion because the coverage is accumulated near the trajectory. The use of the Gaussian function in the HEDAC control has been proven in numerical experiments in two-dimensional coverage motion control applications \cite{ivic2016ergodicity, ivic2020motion, ivic2022constrained}. This  instantaneous action is defined as:
\begin{equation}
	\phi_\sigma(r) = \frac{\Phi}{\left(\sigma \sqrt{2\pi}\right)^3} \exp\left(-\frac{r^2}{2\sigma^2}\right) 
\end{equation}
where $\Phi$ is the coverage action intensity, $\sigma$ is the standard deviation or the scope of action function $\phi$, and $r_i(t) = \|\mathbf{x} - \mathbf{y}_i(t)\|$ is the distance from the location in the domain $\mathbf{x}$ to the agent's position $\mathbf{y}$. The intensity can be interpreted via a property of the instantaneous action function: $\int_{\mathbb{R}^3} \phi \diff\mathbf{x} = \Phi$.
The coverage field $\rho(\mathbf{x})$ is defined as a convolution integral:
\begin{equation}
	\rho(\mathbf{x},t) = \sum_{i=1}^n \int_0^t \phi_\sigma(r_i (\tau )) \diff \tau.
	\label{eq:coverage}
\end{equation}
Note that, due to the simplicity and practicality, the same action function $\phi_\sigma$ is utilized for all agents. Though, it has been shown in \cite{ivic2020motion} that HEDAC can be utilized for governing multiple heterogeneous agents in both sensing (action) and motion characteristics.

The objective of the space coverage algorithm is to produce agents' trajectories that explore the domain according to a given density. Furthermore, we want to enable a continuous non-stopping coverage motion whenever the agents' safety constraints are not violated. According to the given target density $\mu_0(\mathbf{x})$, characterizing UAV visiting areas and frequency, one can define an exponential law to determine the remaining density $\mu$ at the time $t$ according to the achieved coverage $\rho$:
\begin{equation}
	\mu(\mathbf{x},t) = \mu_0(\mathbf{x})\cdot \exp\left(-\rho(\mathbf{x,t})\right).
	\label{eq:mu}
\end{equation}
The target density $\mu_0$ is normalized to satisfy $\int_\Omega\mu_0(\mathbf{x})\diff\mathbf{x}=1$.
This formulation is analogously used in \cite{ivic2020motion} for defining the (undetected) target probability in multi-agent search motion control with uncertain target detection.

The objective of the ergodic exploration is to evenly fill the space with trajectories or to fill it according to a given density. In the proposed formulation \eqref{eq:mu}, minimizing $\mu$ does not lead to equalization of the target density $\mu_0$ and coverage $\rho$. However, it is trivial to show that minimizing $\mu$ results with $\rho(\mathbf{x})=\ln(\mu_0(\mathbf{x}))$, i.e. $\ln(\mu_0)$ can be considered the goal density used in \cite{ivic2016ergodicity, mathew2011metrics}. This approach elegantly solves the problems of normalization and logarithm of zero-valued $\mu_0$, which can arise in the conventional ergodic coverage formulation.
As stated, in order to realize spatial coverage, the motion control needs to direct UAVs to accomplish:
\begin{equation}
	\lim_{t \to \infty} \int_\Omega \mu(\mathbf{x}, t) \diff\mathbf{x} = 0.
	\label{eq:limit_minimization}
\end{equation}
It is suitable to define the measure of spatial density coverage:
\begin{equation}
	\eta_V(t) =  1 - \int_\Omega \mu(\mathbf{x}, t) \diff\mathbf{x},
\end{equation}
which indicates a share of the covered space.

\subsection{Utilizing a potential field for directing UAVs}

In order to minimize \eqref{eq:limit_minimization} in time, one needs to minimize the $\mu$ spatially. HEDAC's main idea is to design a potential field that can facilitate the minimization of $\mu$. The Helmholtz partial differential equation, used for modeling the conductive heat transfer accompanied by convective cooling, is employed for obtaining the potential $\psi$:
\begin{equation}
	k\cdot\Delta\psi(\mathbf{x}. t) - \psi(\mathbf{x}, t)  + \mu(\mathbf{x}, t) = 0
	\label{eq:helmholtz}
\end{equation}
with Neumann boundary condition applied to the entire boundary:
\begin{equation}
	\left.\frac{\partial\psi}{\partial\mathbf{n}}\right|_\Gamma = 0,
	\label{eq:neumann_bc}
\end{equation}
where $\mathbf{n}$ is the boundary outward-pointing normal.
The parameter $k$ represents the coefficient of thermal conductivity and regulates global and local details of the resulting potential field. For the purpose of this search model, we assume this coefficient is dimensionless. Note that more parameters have been used in the previous HEDAC formulations, but the conduction coefficient is dominantly regulating the behavior of $\psi$ hence it is reasonable to consider only this parameter of the HEDAC control.

The potential field $\psi$ calculated by \eqref{eq:helmholtz} is actually a smoothed field $\mu$, which is accomplished due to using the Laplacian operator $\Delta\psi$. This allows us to utilize the gradient of the potential in order to establish the direction to regions of higher potential and implicitly to regions of higher values of density $\mu$. We calculate the desired direction of each UAV motion as a unit gradient of the potential $\psi$:
\begin{equation}
	\mathbf{u}_i(t) = \frac{\nabla\psi(\mathbf{y}_i, t)}{\|\nabla\psi(\mathbf{y}_i, t)\|}, \qquad i=1,\ldots, n.
\end{equation}

Finally, we define a motion model to close the control's feedback loop. 
The motion, and consequentially the trajectories, are an outcome of the 1-st order control ${d \mathbf{y}_i}/{d t} = v_i\cdot\mathbf{u}_i(t)$, where the position is directly changed by the control-appointed direction. The complete motion model, including the collision avoidance mechanism, is defined in subsection \ref{subsec:uavtraj}. 
This motion model is often called kinematic as it neglects the mass and inertial effects of the UAV motion. However, due to advanced low-level multi-rotor control, in practice, modern multi-rotors can make a turn in any direction almost instantaneously. Therefore, it is justified to use such a model for multi-rotor UAV inspection simulations. 

Due to the Neumann boundary condition, the gradient of the potential $\nabla\psi$ inherently prevents an agent from approaching the domain boundaries. However, there is no assurance that the agent will not collide with the boundary due to a chosen combination of numerical parameters for solving the partial differential equation \eqref{eq:helmholtz} (such as numerical grid density), UAV motion time step, and UAV properties (such as velocity $v$ and action range $\sigma$). Furthermore, a possible collision between two or more UAVs also needs to be taken into account when directing multiple UAVs.
\subsection{Collision avoidance}
Several different approaches for boundary and inter-UAV collision avoidance are successfully used in two-dimensional HEDAC applications \cite{ivic2016ergodicity,ivic2019autonomous,ivic2022constrained}. Compared to collision avoidance in two dimensions, the three-dimensional case has more degrees of freedom hence the probability of random collisions between UAVs is greatly reduced. However, they do eventually occur if no collision avoidance procedure is implemented. 
We have implemented a robust and computationally inexpensive collision avoidance mechanism that prevents collisions of agents with other agents and with the domain boundaries.

The possibility of collisions is checked at every time step for $i$-th agent by measuring the distance $d(i,j)=|\mathbf{y}_j-\mathbf{y}_i|$, $i\neq j$ to all other agents indexed with $j$, where $d(i,i)=\infty$, and the minimal distance from the agent to the domain boundary $d_b(\mathbf{y}_i)$. If the minimum of all measured distances $d(\mathbf{y}_i)=min(d_b(\mathbf{y}_i),d(i,1),\ldots,d(i,n))$ is below initially prescribed threshold $2\epsilon$,  the collision avoidance must intervene and correct the agent's direction vector. To ensure that the outlines of a UAV do not collide with obstacles, the safety distance must be larger than the maximum dimension from the center of the UAV $M$ plus the distance that the UAV can reach in a time step. Safety distance parameter  $\epsilon>v_i \Delta t + M$ must be large enough to enable safe and continuous maneuvers when avoiding collisions and low enough to not interfere significantly with the primary goal.  If the agent's distance to obstacles is below the safety distance, it must not be allowed to reduce this distance in the next step. 

Let $I_i=\{ j|d(i,j)<2\epsilon, i\neq j \}$ 
be a set of indices of agents for which the $i$-th agent is in danger of collision and suppose that $d_b(\mathbf{y}_i)<2\epsilon$ i.e. the $i$-th agent is also to close too the domain boundary. We would like to find a direction vector $\mathbf{w}_i$ for the $i$-th agent such that it closes an obtuse angle with vectors $\mathbf{e}_{i,j}=\mathbf{y}_j-\mathbf{y}_i,j\in I_i$ and vector $\nabla d_b(\mathbf{y}_i)$. Furthermore, to make this vector unique, we can pick the one which is closest to the vector $\mathbf{c}=\sum_{j\in I} \mathbf{e}_{i,j}/\|\mathbf{e}_{i,j}\|+\nabla d_b(\mathbf{y}_i)/\|\nabla d_b(\mathbf{y}_i)\|+\mathbf{u}_i$.

If we are able to solve this problem, then this direction will lead UAVs away from the boundary and neighboring agents in the next step. This problem can be formulated as an optimization problem: 

\begin{align}
	& \text{maximize}_{\mathbf{w}_i}   && \mathbf{c}^T \mathbf{w}_i\\
	& \text{subject to} && -\mathbf{A} \mathbf{w}_i \leq \mathbf{0} \\
	&  && ||\mathbf{w}_i|| = 1,
	\label{eq: nonlinearC}
\end{align}
where 
\begin{equation*}
	\mathbf{A}=
	\begin{bmatrix}
		\bar{\mathbf{A}}\\
		\nabla d_b(\mathbf{y}_i)\\ 
		\mathbf{u}_i
	\end{bmatrix}
\end{equation*}
and $\bar{\mathbf{A}}=[\mathbf{e}_{i,j}],j\in I_i$.

\begin{figure}[!htb]
	\centering
	\includegraphics[width=\linewidth]{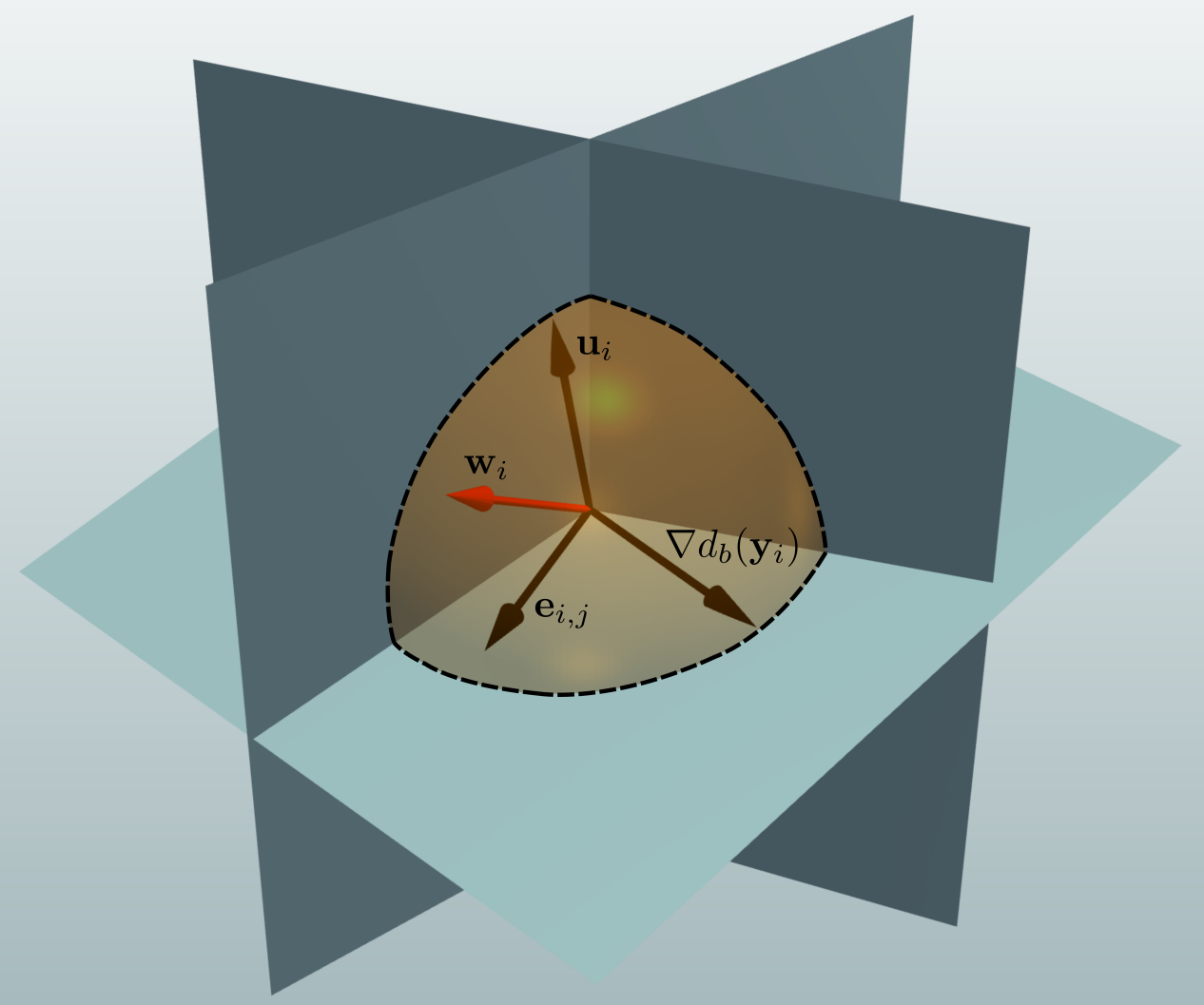}
	\caption{Feasible region for the optimal collision avoidance vector marked with a dashed sphere triangle. Solution $\mathbf{w}_i$ of the optimization problem with nonlinear constraint $||\mathbf{w}_i|| = 1$ and linear constraints defined by normal vectors $\mathbf{u}_i$, $\mathbf{e}_{i,j}$ and $\nabla d_b(\mathbf{y}_i)$.}
	\label{fig:minimizaiton}
\end{figure}

The solution of this linear objective function with nonlinear constraints can be solved efficiently using a nonlinear solver. Figure \ref{fig:minimizaiton} shows the feasible space for a matrix $A$ with three rows where each row corresponds to a depicted black vector pointing towards the feasible region. The solution is a unit vector $\mathbf{w}_i$ depicted in red which closes an acute angle with highlighted vectors. 

It can be shown that the distance from the boundary $d_b\in C^2$, therefore, $\nabla d_b(\mathbf{y}_i)$ always exists and points towards the boundary. We should underline that the matrix $\mathbf{A}$ has at least one row and $\nabla d_b(\mathbf{y}_i)$ is not included if $d_b(\mathbf{y}_i)\geq \epsilon$. It is possible that the optimization problem does not have a solution because the feasible set is empty. In that case, the matrix $\mathbf{A}$ must have at least 6 rows (conditions) which is highly unlikely in practical situations. In the two-dimensional case, it would take at least 4 rows for this to occur.
In the unlikely case that the optimization problem is not solved, the UAV should stand still for this time step and wait for the neighboring agents to clear its path. This maneuver is feasible for a multi-rotor UAV.

An initial approximate solution to the optimization problem is obtained using an approximate linear program that finds a good initial solution candidate in the feasible space. The initial approximate solution is used by a Trust-region optimization implemented in Scipy 1.7.3 \cite{2020SciPy-NMeth}  which is based on algorithms found in \cite{Conn2000}.
The solution to this optimization problem is computationally inexpensive and is invoked sparsely therefore it does not add significantly to the overall complexity of the algorithm.

\subsection{Defining UAVs' trajectories}
\label{subsec:uavtraj}
The motion of each agent in the fleet of UAVs is defined as: 
\begin{equation}
	\frac{d \mathbf{y}_i}{d t} = 
	\begin{cases}
		v_i\cdot\mathbf{w}_i(t) & \text{ if } d(\mathbf{y}_i) < \epsilon,\\ 
		v_i\cdot \left(\left(2-\frac{d(\mathbf{y}_i)}{\epsilon}\right)\mathbf{w}_i(t)+ \left(\frac{d(\mathbf{y}_i)}{\epsilon}-1\right) \mathbf{u}_i(t)\right)& \text{ if } \epsilon \leq d(\mathbf{y}_i) < 2\epsilon,\\ 
		v_i\cdot\mathbf{u}_i(t) & \text{otherwise} ,\\
	\end{cases}
\end{equation}

where $v$ is a velocity magnitude equal for all UAVs, and $\mathbf{w}_i$ is the direction obtained as a solution of the nonlinear optimization problem.

The collision avoidance can cause jitters in the trajectories because it is not applied smoothly. To fix this problem and smooth out trajectories, the correction of the direction vector is applied gradually when the minimal distance satisfies $d(\mathbf{y}_i)\in [\epsilon,2\epsilon)$. In the unlikely case that the agent gets too close to an obstacle, the collision-avoidance direction vector $\mathbf{w}_i$ is applied to move the agent toward a collision-safe part of the domain.

\subsection{Numerical implementation using finite element method}
The numerical solution of equation \eqref{eq:helmholtz} together with boundary condition \eqref{eq:neumann_bc} can be complicated to solve, especially if the geometry of the boundary is complex. In practical applications of this model, we use a connected 3D domain with holes that can have complex geometry. Furthermore, for practical applications, we need a relatively large number of numerical grid points. FEM provides a simple implementation of boundary conditions, a fast solver, and a straightforward interpolation of results on any given grid.

The weak formulation of the presented problem is obtained by multiplying the equation \eqref{eq:helmholtz} by a smooth test function $v\in H^1(\Omega),\Omega\subset\mathbb{R}^3$ and integrating over the domain. Using integration by parts, the following equation is obtained:
\begin{eqnarray*}
	-k \int_{\Omega} \nabla \psi (\mathbf{x}, t)\nabla v (\mathbf{x})\diff \Omega - \int_{\Omega} \psi (\mathbf{x}, t) v (\mathbf{x})\diff \Omega \\
	+ \int_{\Omega} \mu (\mathbf{x}, t)v (\mathbf{x})\diff \Omega  +\sum_{j=0}^{n_0}\oint_{\Gamma_j} (\nabla \psi (\mathbf{x}, t). \mathbf{n})v(\mathbf{x})\diff \Gamma= 0.
\end{eqnarray*} 
The final form of the weak formulation is obtained after we apply the Neumann boundary for all parts of the boundary $\Gamma_j\subset\Gamma$:
\begin{equation}
	\int_{\Omega} k\nabla \psi (\mathbf{x}, t)\nabla v (\mathbf{x}) + \psi (\mathbf{x}, t) v (\mathbf{x}) -  \mu (\mathbf{x}, t)v (\mathbf{x})\diff \Omega = 0.
	\label{eq:hedac_week}
\end{equation} 

We use quadratic polynomials for the space of test functions and our representation of the solution.
Let $M=\left\{T_{1}, \ldots, T_{N}\right\}$ be a partition of $\Omega$ into $N$ uniform non-overlapping triangles. The triangles with the geometry are described by the classical 3-node interpolation functions. The scalar field of the unknown variable and test functions over each $n$-node element $T_{i}$ is approximated by
$$\psi_{h}=\sum_{j=1}^{n} N_{j} u_{j}$$
where $N_{j}$ stands for the Lagrangian polynomial interpolation functions. 

Because the domain and the triangulation do not change during the entire calculation, the linear system coefficient matrix is sparse and constant hence the solutions are obtained very efficiently. The interpolation of different scalar fields and calculation of the gradient $\nabla \psi$ can be directly obtained using a finite element representation \cite{NGSolve}. 
All results are obtained with the latest version of the finite element software NGSolve.

\subsection{3D coverage example}
The proposed algorithm can be used to cover a three-dimensional space with a given density of trajectories. To demonstrate this application we use a simple unit cube domain and $100$ agents initially distributed randomly in the center of the cube. The density of the trajectories must uniformly fill the lower part of the cube where $0.05<x,y<0.95$ and $0.05 < z < 0.45$. This test demonstrates that this algorithm can avoid collisions inside a crowded domain without sacrificing the goal of space coverage. The safety distance is set to $0.025$ m and the velocity is a constant of $0.1~\text{m/s}$. Total simulation time is $2000~\text{s}$ with the time step $\Delta t=1~\text{s}$ and the conduction parameter $k = 0.1$. 

\begin{figure}[!htb]
	\centering
	\raggedright \footnotesize ~(A)
	\includegraphics[width=\linewidth]{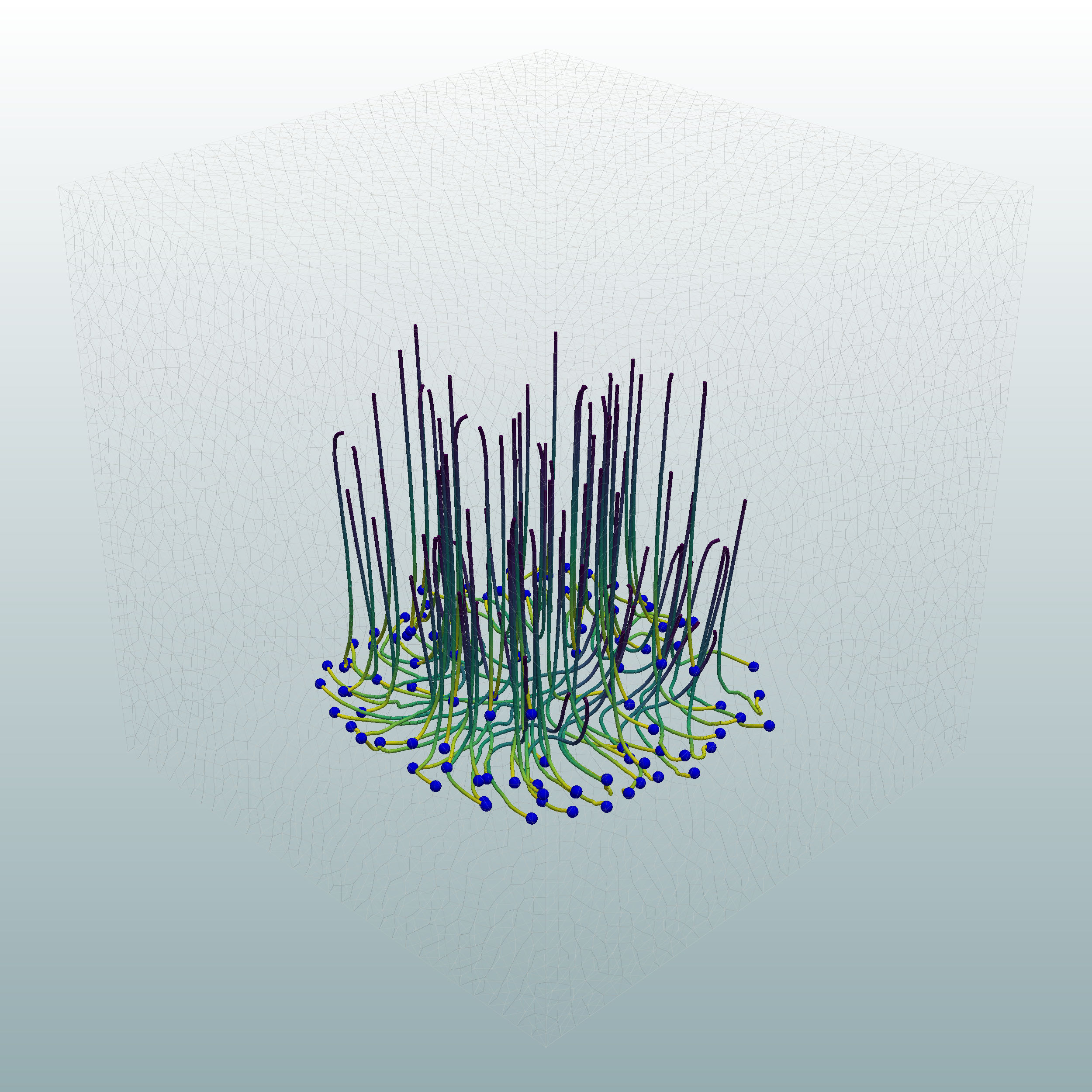}
	\raggedright \footnotesize ~(B)
	\includegraphics[width=\linewidth]{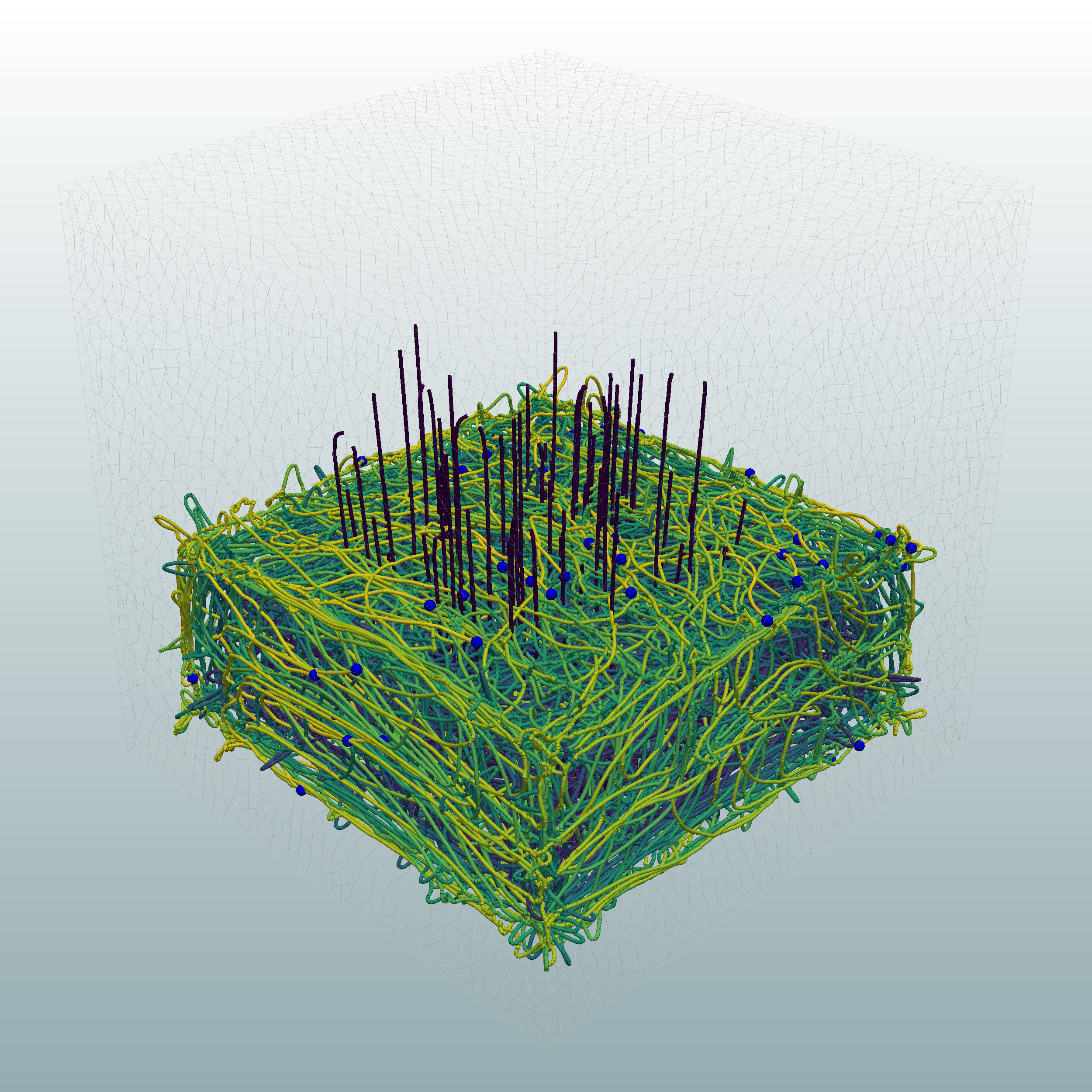}
	\caption{Trajectories of 100 agents after $50$ (A) and $1000$ seconds (B) produced with HEDAC method for uniform coverage of lower half of a unit cube.}
	\label{fig:cube200}
\end{figure}

Figure \ref{fig:cube200} shows the final result after $50$ s and $1000$ s. One can observe that the trajectories form a space-filling curve and the agents keep inside the goal area except for a short initial part of the simulation. This domain is crowded with agents but the trajectories are quite smooth and no deadlocks or jitters in the produced trajectories can be observed. 

\begin{figure}[!htb]
	\centering	
	\includegraphics[width=\linewidth]{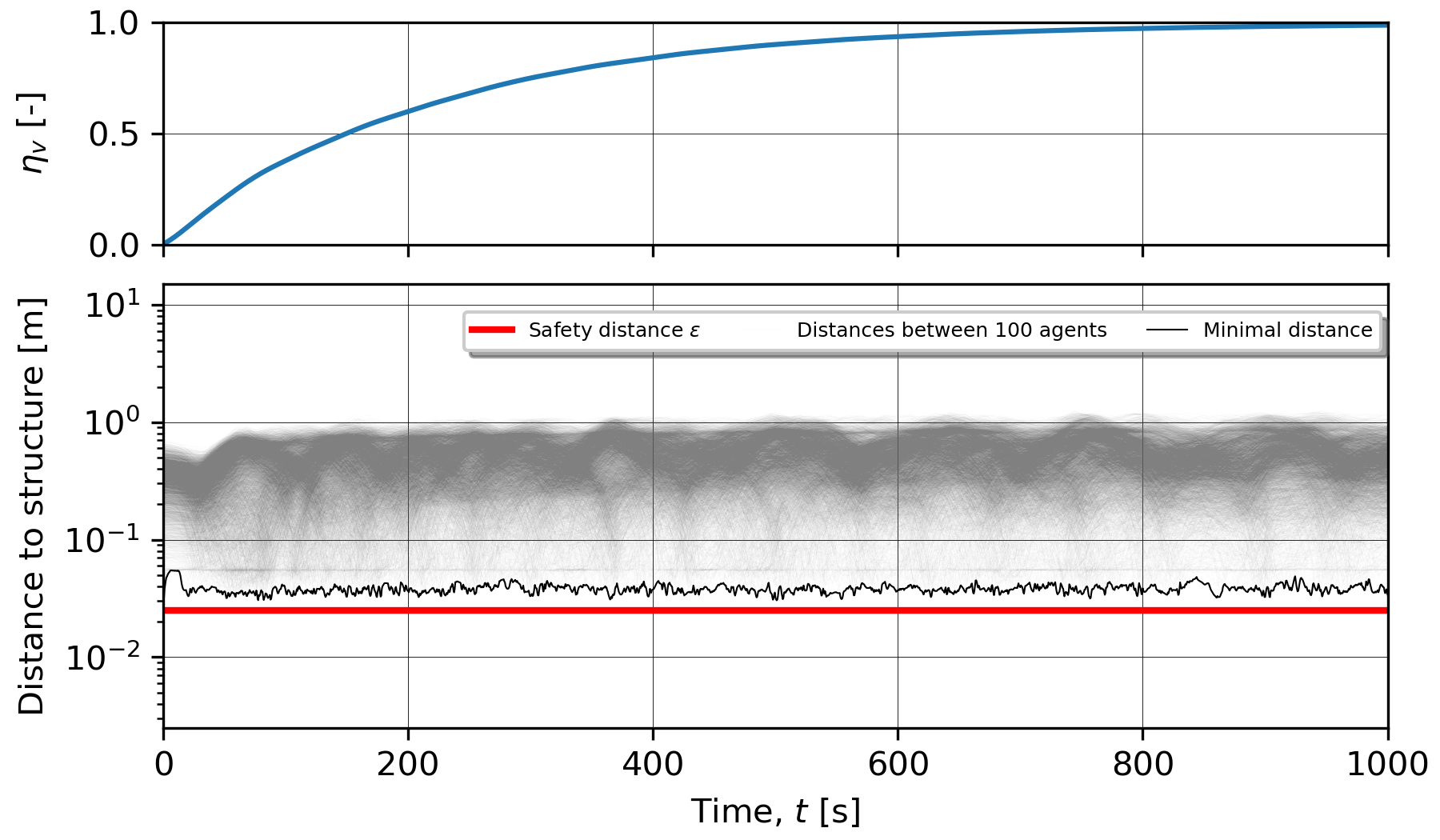}	
	\caption{Convergence of the spatial density coverage measure $\eta_V$ and computed UAV distances during the inspection of a simple unit cube with $100$ agents.}
	\label{fig:case_1_unit_cube_convergence}
\end{figure}

Figure \ref{fig:case_1_unit_cube_convergence} clearly shows that the given domain is explored efficiently and the collision avoidance is working properly in this crowded domain. The agents safely navigate the domain and avoid other agents in the process. We must emphasize that the agents, in this case, are more in danger of colliding with each other than the domain boundary.
This is a simple artificial test example, but the results demonstrate the algorithm is ready to be tested in a more realistic setting with complex geometries.

\section{UAV visual inspection applications}
\label{sec:four}

The main task in the visual inspection is to capture all surfaces of the structure using a camera mounted on the UAV. This requires control of the direction of the camera depending on the position of the drone relative to the observed structure. We propose a rather simple idea that considers the control of camera orientation, along with obtained UAV trajectory, by directing the camera view to the nearest point of the observed structure. 

In order to demonstrate the applicability of the proposed algorithm to infrastructure inspection, we conduct inspection simulations for three different test cases with varying complexity.

\subsection{Area of interest and camera control for 3D structure inspection}
We propose a camera direction control that can provide relatively good results if trajectories obtained by HEDAC suitably explore the region around the structure.
Since the inspection is acting on faces and surfaces of the investigated structure, the scope of required trajectories can not be correctly transferred to the three-dimensional domain $\Omega$. However, the proposed approximation is satisfactory as it is shown in the UAV inspection scenarios presented in this section and further investigated in Section \ref{sec:five}.

Similar to the $d_b$ field, we use $d_s(\mathbf{x})$ as the distance of $\mathbf{x}$ to the structure surface. Now, we can define a suitable region of interest around the observed structure using a three-dimensional Gaussian function:
\begin{equation}
	\mu_0(\mathbf{x}) = \exp\left(-\dfrac{\left(d_s(\mathbf{x}) - {d_m}\right)^2}{2 \cdot {d_\sigma}^2}\right),
	\label{eq:mu0_for_inspection}
\end{equation}
where ${d_m}$ is the goal inspection distance at which the Gaussian function is centered and ${d_\sigma}$  is the standard deviation, i.e. the broadness of the field $\mu_0$ encompassing the structure. Using \eqref{eq:mu0_for_inspection} allows us to construct a continuous non-negative field around the entire structure, with a peak exactly at the distance ${d_m}$ from the surfaces to be inspected and gradually weakening both towards and away from the structure. Note that scaling constants are omitted from the Gaussian function in \eqref{eq:mu0_for_inspection} since $\mu_0$ is subject to normalization that produces compliant scaling constants. 

Finally, we can easily find a camera orientation $\mathbf{z}_i$ for each UAV, as the direction towards the nearest point on the structure, by employing the gradient of the field $d_s$:
\begin{equation}
	\mathbf{z}_i(\mathbf{y}_i) = \frac{\nabla d_s(\mathbf{y}_i)}{\| \nabla d_s(\mathbf{y}_i) \|}.
	\label{eq:camera_orientation}
\end{equation}
Note that only camera orientation is subjected to the proposed control and provided by the unit vector $\mathbf{z}$. This approach does not acknowledge the field of view (FOV), camera focus, zoom, or other details regarding photographic equipment for visual inspection. Although not considered in the trajectory planning algorithm, a simplified field of view has been utilized for surface coverage assessment as described in the next subsection.

\subsection{Surface coverage assessment}

\begin{figure}[!htb]
	\centering	
	\includegraphics[width=\linewidth]{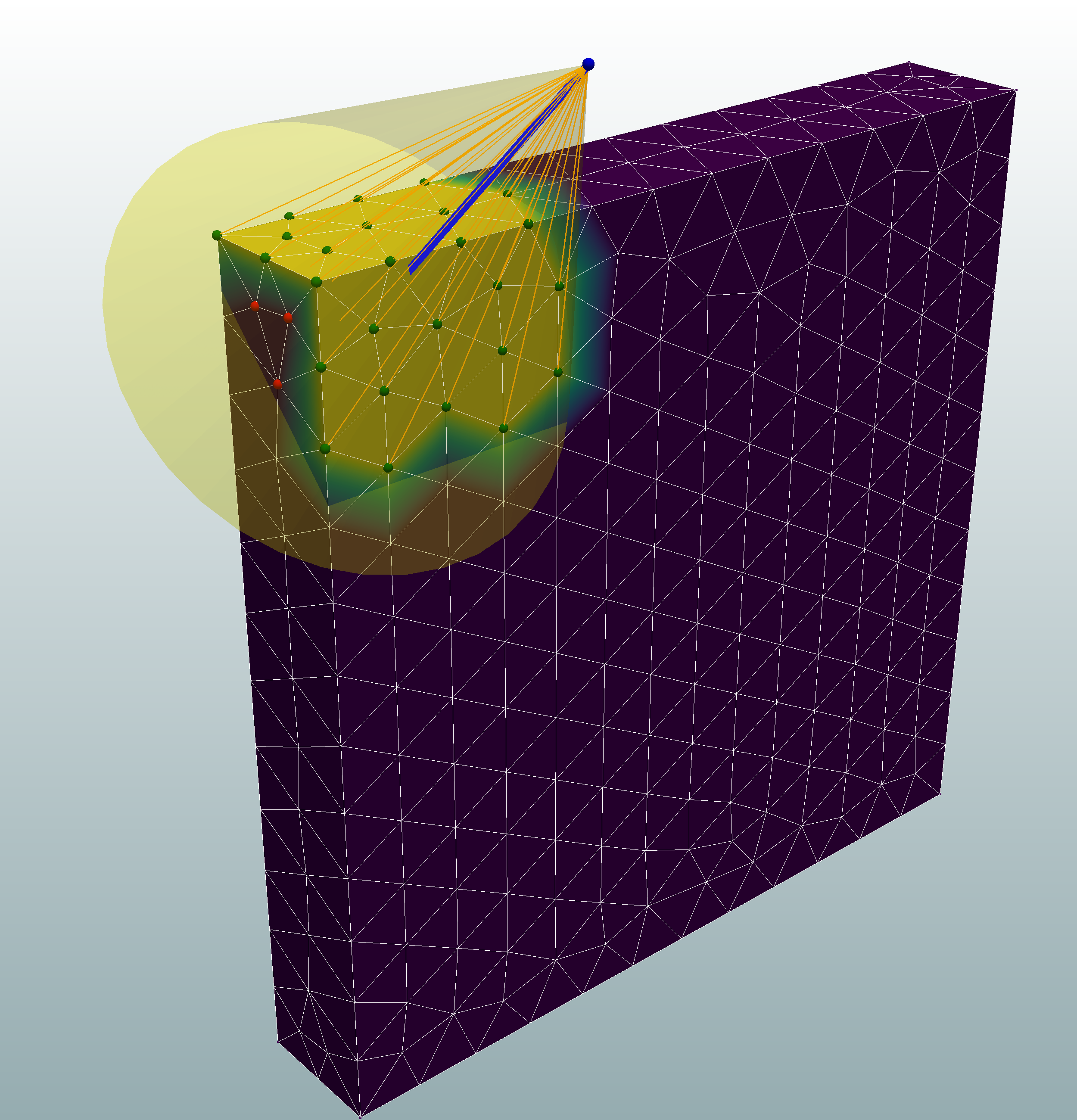}	
	\caption{Example of surface inspection evaluation via calculation of observable nodes (marked as green dots) in a conical field of view (semi-transparent yellow cone). Note that not all candidate nodes (nodes inside of the cone) are marked as observed. Nodes marked as red dots did not pass a "direct line of sight" test, i.e. their orange ray intersects the surface before landing into the candidate node. The blue dot represents the UAV/camera location while the thick blue line represents the camera orientation $\mathbf{z}$ i.e. the cone axis.}
	\label{fig:surface_coverage_calculation}
\end{figure}

In order to validate the correctness of solving the surface coverage task via ergodic exploration of the spatial field, we implement a simulation of the inspection camera view (Figure~\ref{fig:surface_coverage_calculation}). The camera's FOV is modeled as a cone whose top is positioned at the UAV's center $\mathbf{y}_i$. The cone orientation corresponds to camera orientation $\mathbf{z}_i$ obtained by \eqref{eq:camera_orientation}. The height of the cone $C_H$ regulates the maximum acceptable distance at which inspection images are taken while $C_R$ is the radius of the FOV at that distance. We only consider a circular area of the image, in contrast to the rectangular image obtainable by the camera. This can be justified by considering only inspection details at the center of the recorded image.

At each time step, all nodes on the structure's surface are tested if they are inside the cone. The candidate nodes (the ones inside the cone) are further tested to check if they are in the camera's line of sight. This is performed by a simple ray tracing technique: if there is no intersection between a ray shot from the camera and the structure's surface, then a candidate node is directly observed. The number of observations is tracked for all nodes during the entire inspection simulation. For relatively uniform meshes, which are used in the presented test cases, the ratio between the number of inspected nodes and the total number of nodes on the structure's surface represents a share of the inspected surface $\eta_A$. Although simplified to ensure computational feasibility, this coverage assesment model provides a qualitatively satisfactory estimate of surface inspection.  

%\section{UAV inspection simulation results}
%\label{sec:five}

\subsection{Portal test case}

\begin{figure}[!htb]
	\centering
	\includegraphics[width=\linewidth]{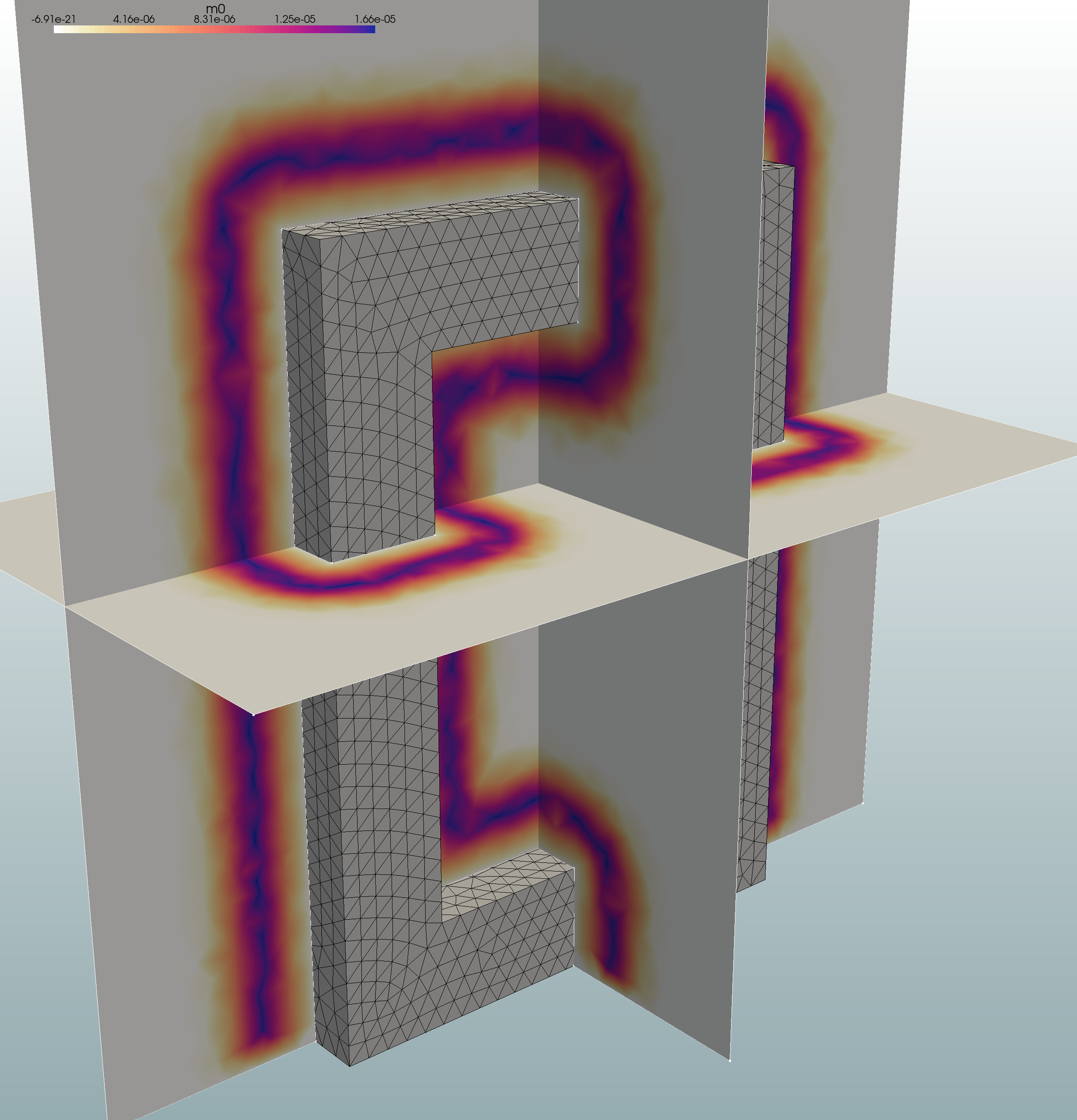}
	\caption{Target density field $\mu_0$, of thickness $d_\sigma=2~\text{m}$, encompasses the portal structure at the distance $d_m=5~\text{m}$.}
	\label{fig:case_2_portal_mu0}
\end{figure}

The portal test case is a scenario for visual inspection of a relatively simple synthetic three-dimensional structure. It can be described as a flattened box ($50 \times 70 \times 10$ m) with a rectangular hole in it ($30 \times 50$ m), and it is designed to provide a simple shape while requiring relatively complex maneuvers to achieve structure inspection. 
The target density field $\mu_0$ (Figure~\ref{fig:case_2_portal_mu0}) is defined using the goal inspection distance $d_m=5~\text{m}$ and the broadness $d_\sigma=2~\text{m}$. Other info and parameters used for the portal case simulation are presented in Table~\ref{tab:case_2_portal_parameters}.

\begin{table}[!htb]
	\small
	\centering
	\caption{Domain, structure and numerical mesh info, inspection task parameters, UAV parameters, and HEDAC parameters used in the portal test case.}
	\begin{tabularx}{\linewidth}{Xrl}
		Parameter & Value & Unit\\
		\hline
		Structure length & 10 & m \\
		Structure width & 50 & m \\
		Structure height & 70 & m \\
		Domain length & 50 & m \\
		Domain width & 90 & m \\
		Domain height & 90 & m \\
		\hline
		Number of domain mesh nodes & 34 161 \\
		Number of domain mesh elements & 167 495 \\
		Number of structure mesh surface nodes & 2 315 \\
		Number of structure mesh surface faces & 4 630 \\
		\hline
		Inspection distance $d_m$ & 5 & m \\
		Inspection distance broadness $d_\sigma$ & 2 & m \\
		\hline
		\hspace{0pt}FOV cone height $C_H$ & 8 & m \\
		\hspace{0pt}FOV cone diameter $C_D$ & 10 & m \\
		\hline
		Number of UAVs & 3 \\
		UAV velocity $v$ & 0.5 & m/s \\
		Safety distance $\epsilon$ & 1 & m \\	
		Coverage action intensity $\Phi$ & 200 \\
		Coverage action range $\sigma$ & 5 & m \\	
		\hline
		HEDAC conduction coefficient $k$ & 12 \\
		Inspection duration & 1 500 & s \\
		Path planning time step $\Delta t$ & 1 & s \\
		\hline
	\end{tabularx}
	\label{tab:case_2_portal_parameters}
\end{table}

\begin{figure}[!htb]
	\centering	
	\raggedright \footnotesize ~(A)
	\includegraphics[width=\linewidth]{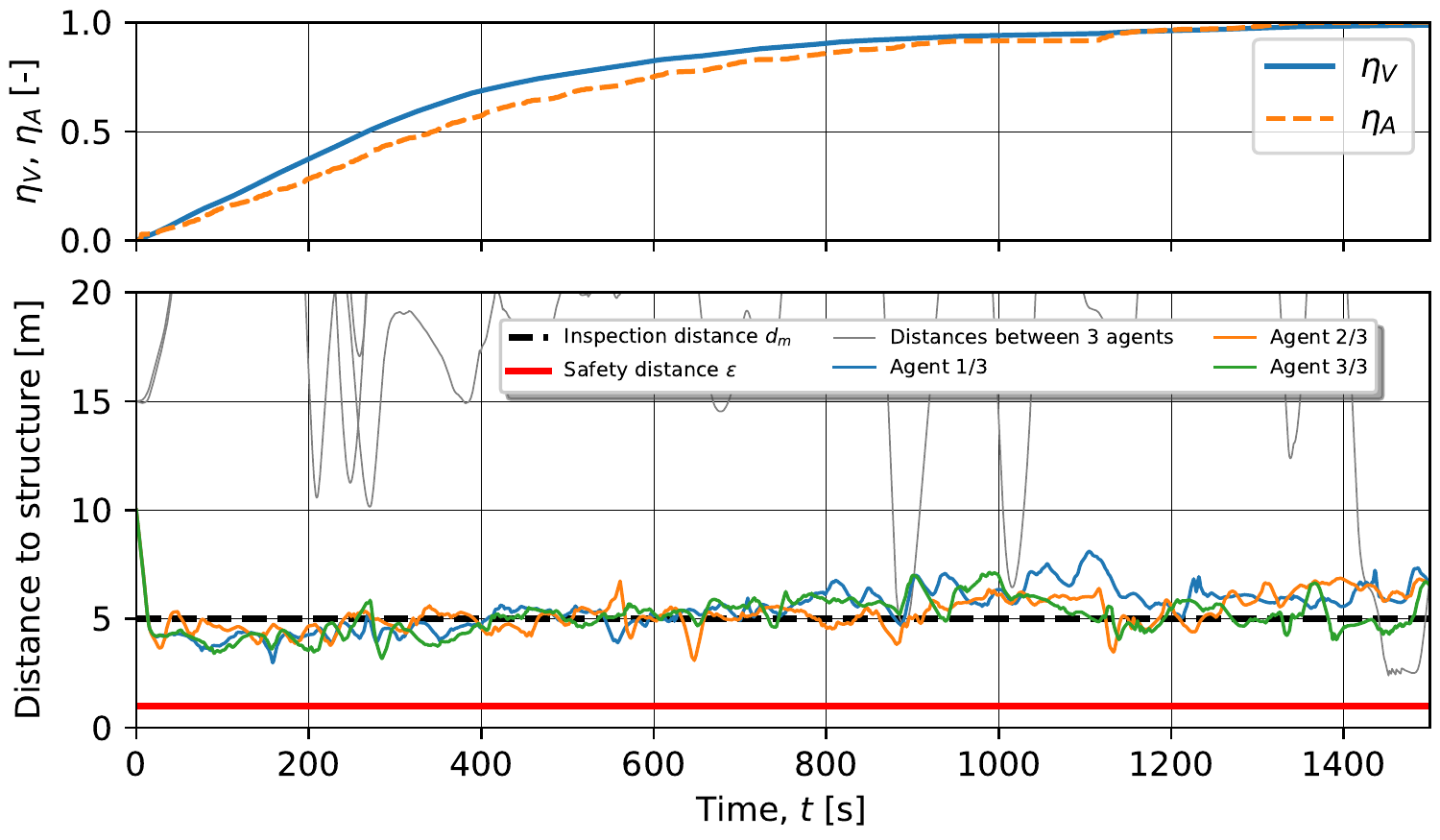}	
	\raggedright \footnotesize ~(B)
	\includegraphics[width=\linewidth]{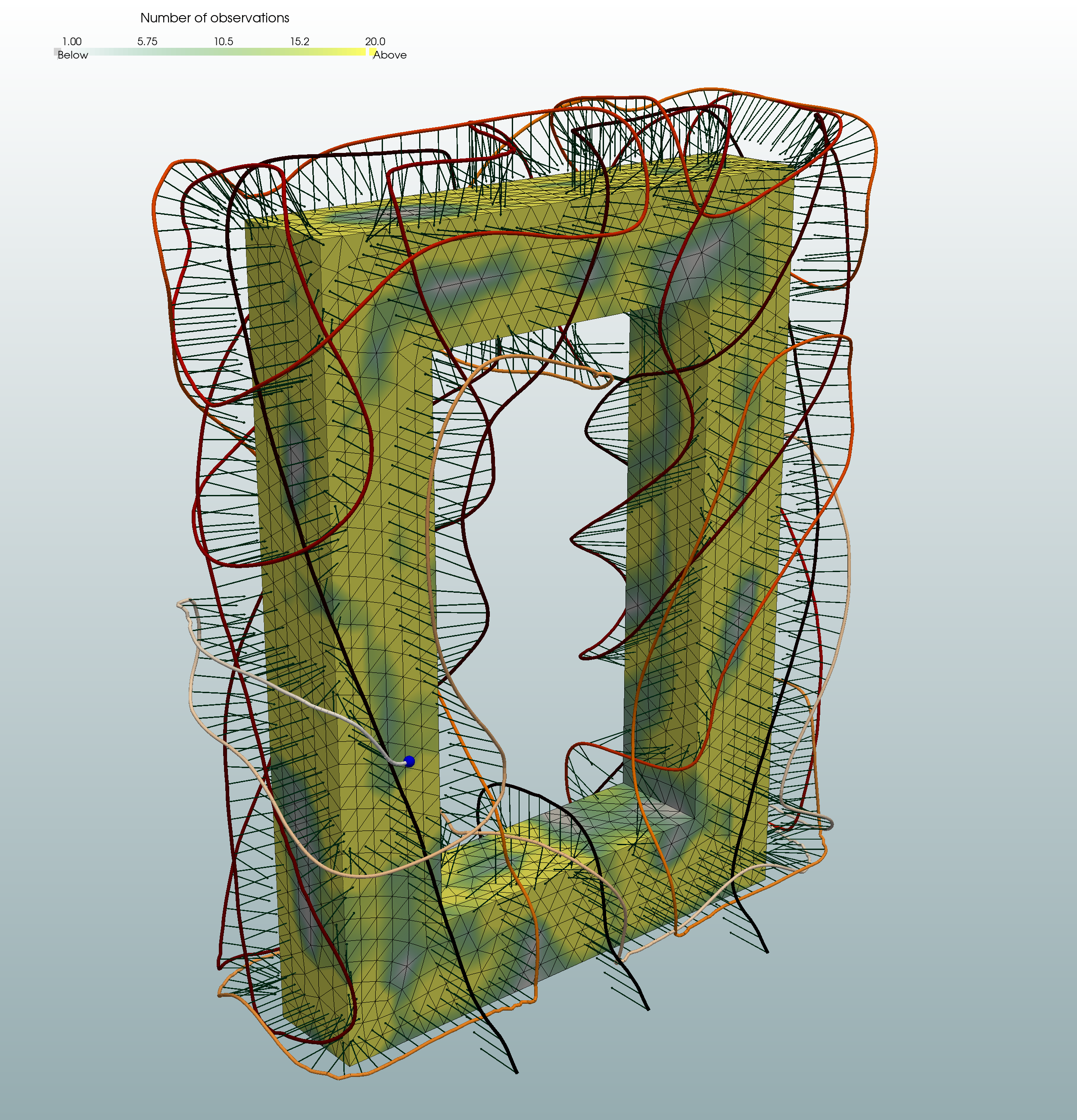}
	\caption{(A) Convergence of the spatial density coverage measure $\eta_{V}$, surface coverage $\eta_A$, and computed UAV distances during the inspection. (B) shows trajectories of three UAVs performing inspection of the portal structure after 1 500 s. Camera orientations along trajectories are displayed as magenta arrows of length equal to the inspection distance $d_m$. Surface inspection is indicated by the number of surface observations visualized as green-to-yellow shading of the portal structure surface.}
	\label{fig:case_2_portal_inspection}
\end{figure}

The proposed trajectory planning algorithm applied for the visual inspection of the portal structure is run with 3 UAVs in the duration of 1 500 s. For the given $\mu_0$, the algorithm achieves convergence of the spatial density and the surface inspection coverage as shown in $\eta_V$ and $\eta_A$ plots in Figure~\ref{fig:case_2_portal_inspection}(A). 

Plotted distance lines indicate that generated UAV paths around the structure are, on average, at distance $d_m=5~\text{m}$ from the structure. 
However, it can be observed that distances slightly increase, on average, during the inspection. 
The target density $\mu_0$ is computed to form a region at the offset from the structure's surface which is mostly convex. The resulting potential field $\psi$ is positioned slightly closer to the structure than $\mu_0$ on the convex part of the structure's surface, due to the smoothing effect of the Laplacian operator in the Helmholtz equation. Consequently, the outcome is UAV trajectories that are closer to the structure than appointed $d_m$. After initial rough passes around the structure, which are relatively close to it, the remaining density $\mu$ is now more accumulated at a distance greater than $d_m$ from the inspected surfaces. Thus, in the second part of the inspection operation, UAV trajectories are generated at a greater distance from the structure. 
The distances of the produced trajectories are safely kept above the minimal safety distance $\epsilon$ during the entire operation.
Due to the large proportions of structure and domain, distances between UAVs are practically negligible in the context of the minimum allowable spacing constraint.

The produced trajectories $\mathbf{y}_i$ and accompanying camera orientations $\mathbf{z}_i$ for all three UAVs are displayed in Figure~\ref{fig:case_2_portal_inspection}(B). Note that the unit vector of camera orientation is scaled to the length equal to the inspection distance $d_m$. The produced trajectories are fairly smooth and suitably distanced from the structure, complying with the given target density field $\mu_0$ presented in Figure~\ref{fig:case_2_portal_mu0}. 

Based on the camera orientation visualization and the surface observation shading shown in Figure~\ref{fig:case_2_portal_inspection}(B), all surfaces of the inspected structure are observed relatively uniformly which brings the conclusion that planned paths and camera orientations are suitable for inspection application.

\subsection{Wind turbine test case}

Wind turbine inspection with UAV technology has become the norm in the past several years as it provides a lot of benefits over manual or ground-based inspection. The main advantages include human safety and high-quality inspection data with a noticeable reduction in data acquisition time. Furthermore, the given benefits of the UAV-based inspection approach enable a reduction in maintenance costs of wind turbine farms \cite{poleo2021estimating}. Due to the given rationale, a number of previous studies have included a wind turbine test case to assess the quality and efficiency of their respective visual inspection algorithms \cite{stokkeland2015autonomous, mansouri2018cooperative, schafer2016multicopter,bircher2016three}, and hence this paper does the same. 

\begin{table}[!htb]
	\small
	\centering
	\caption{Parameters for the wind turbine inspection scenario containing domain, structure and numerical mesh info, inspection task parameters, UAV parameters, and HEDAC parameters.}
	\begin{tabularx}{\linewidth}{Xrl}
		Parameter & Value & Unit\\
		\hline
		Structure length & 25.6 & m \\
		Structure width & 122.8 & m \\
		Structure height & 203.1 & m \\
		Domain length & 65.9 & m \\
		Domain width & 162.9 & m \\
		Domain height & 223.1 & m \\
		\hline
		Number of domain mesh nodes & 124 979 \\
		Number of domain mesh elements & 711 003 \\
		Number of structure mesh surface nodes & 11 052 \\
		Number of structure mesh surface faces & 22 100 \\
		\hline
		Inspection distance $d_m$ & 6 & m \\
		Inspection distance broadness $d_\sigma$ & 1 & m \\
		\hline
		\hspace{0pt}FOV cone height $C_H$ & 8 & m \\
		\hspace{0pt}FOV cone diameter $C_D$ & 6 & m \\
		\hline
		Number of UAVs & 2 \\
		UAV velocity $v$ & 1.2 & m/s \\
		Safety distance $\epsilon$& 1 & m \\	
		Coverage action intensity $\Phi$ & 400 \\
		Coverage action range $\sigma$ & 3 & m \\	
		\hline
		HEDAC conduction coefficient $k$ & 5 \\
		Inspection duration & 1 200 & s \\
		Trajectory planning time step $\Delta t$ & 0.5 & s \\
		\hline
	\end{tabularx}
	\label{tab:case_3_wind_turbine_parameters}
\end{table}

\begin{figure}[!htb]
	\centering
	\includegraphics[width=\linewidth]{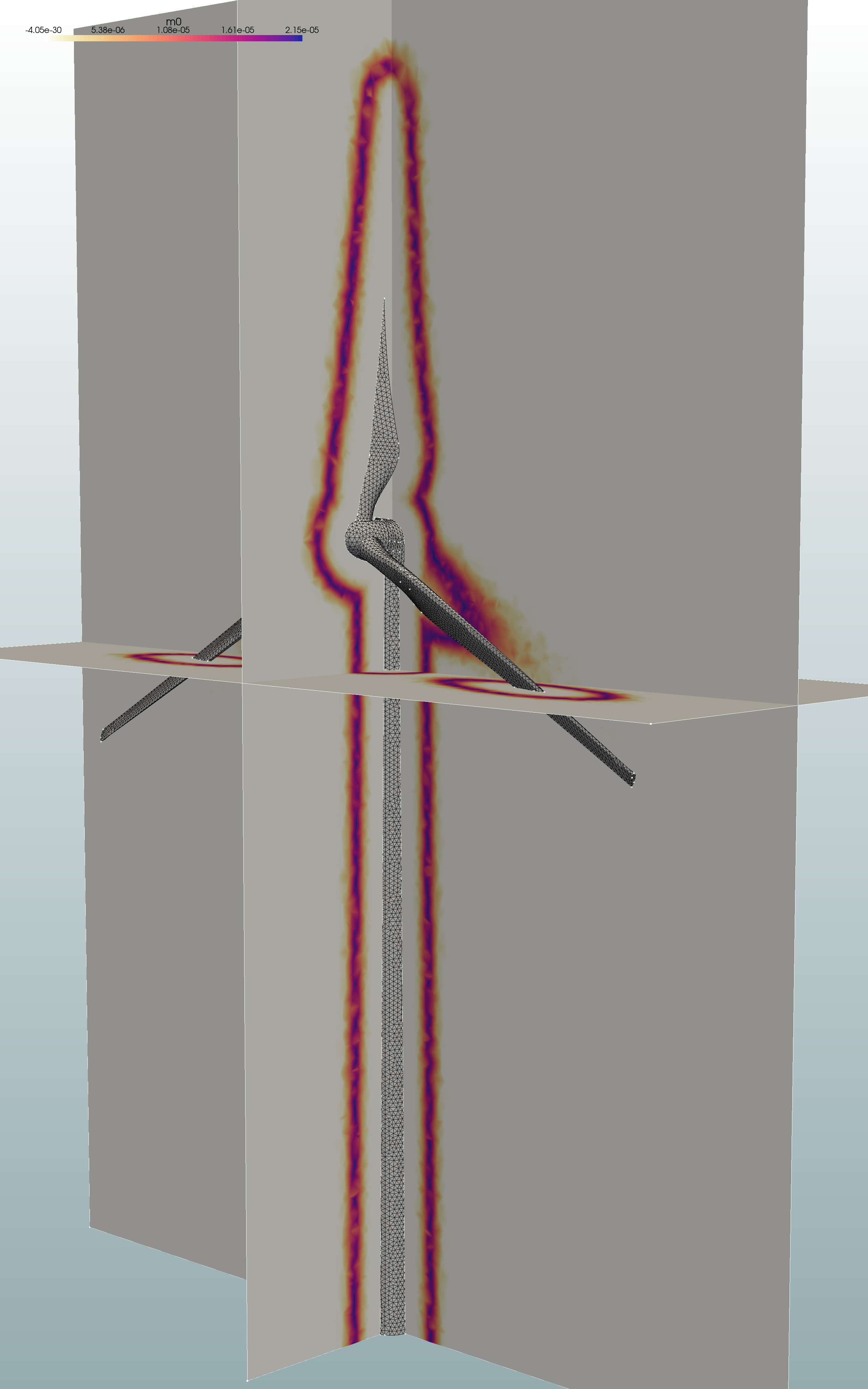}
	\caption{Field $\mu_0$ for the wind turbine test case defines the region aimed to be occupied by UAVs' trajectories along which all turbine structure faces can be inspected.}
	\label{fig:case_3_wind_turbine_mu0}
\end{figure}

For the wind turbine test case, a path planning for the inspection using two UAVs is prepared according to the parameters shown in Table~\ref{tab:case_3_wind_turbine_parameters}. The structure of the wind turbine consists of a 120 m high vertical column, a hub, and three attached blades, each 80 m long. Blades are rotated at a 60$^{\circ}$ angle, forming an "upside-down Y" shape, which is a typical orientation for inspection operations. Based on the shape of the wind turbine, the inspection distance $d_m$ = 6 m and the broadness $d_\sigma$ = 1 m, the target density field $\mu_0$ is computed (Figure~\ref{fig:case_3_wind_turbine_mu0}). Note that the numerical domain used in this test case is significantly larger than the space needed for inspection flight in order to allow for the UAV to directly pass from the tip of one blade to the tip of another blade (or to the root of the column), though this possibility is not utilized in performed trajectory planning computations.

The wind turbine inspection typically takes about 40 minutes during which about 1000 photographs of the turbine surface are recorded. In this inspection scenario, two UAVs are used in order to speed up the inspection operation (the duration is 20 minutes) and to demonstrate the coordination between two UAVs on a relatively slender structure such as a wind turbine.

\begin{figure}[!htb]
	\centering	
	\raggedright \footnotesize ~(A)
	\includegraphics[width=\linewidth]{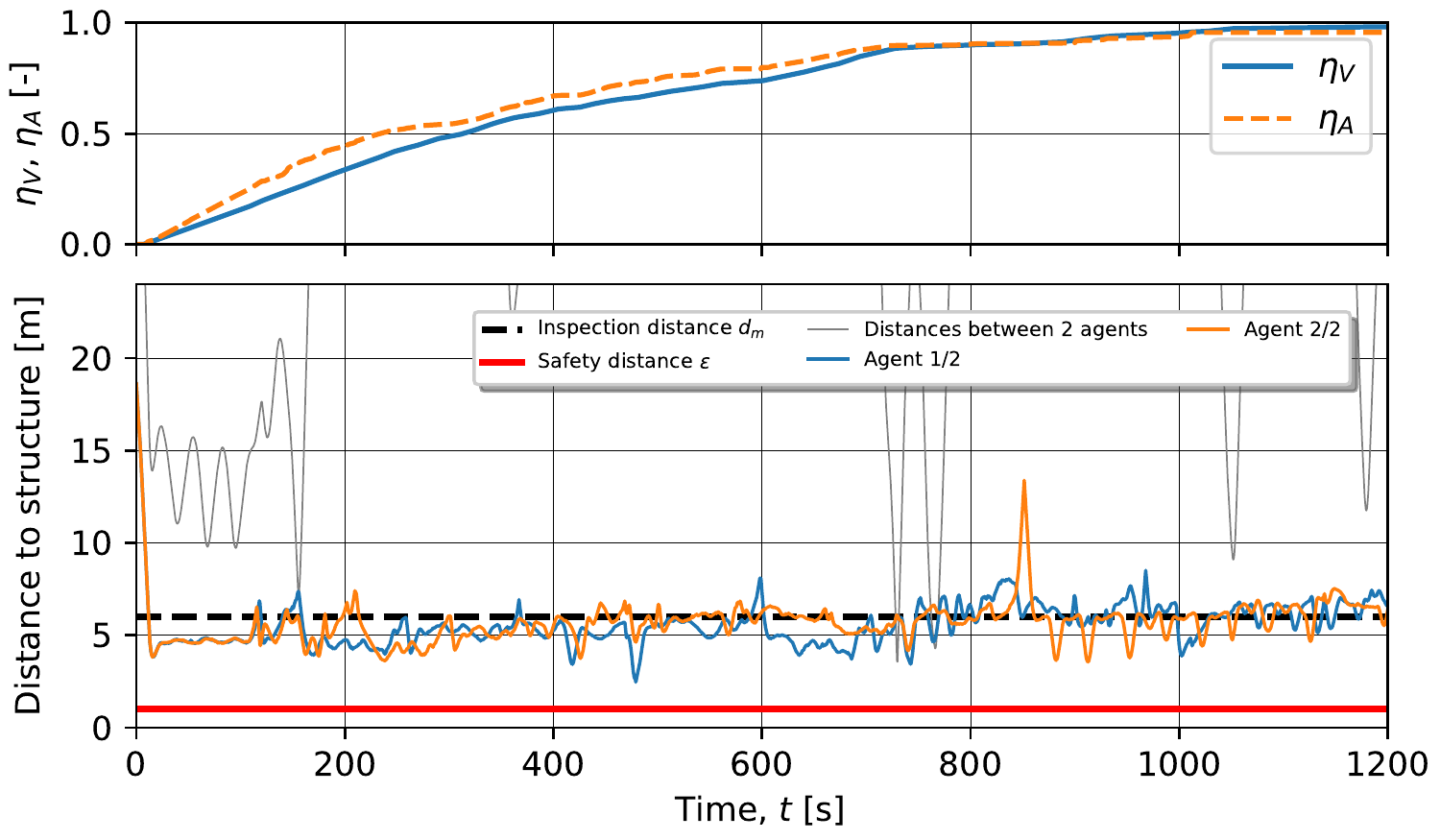}	
	\raggedright \footnotesize ~(B)
	\includegraphics[width=\linewidth]{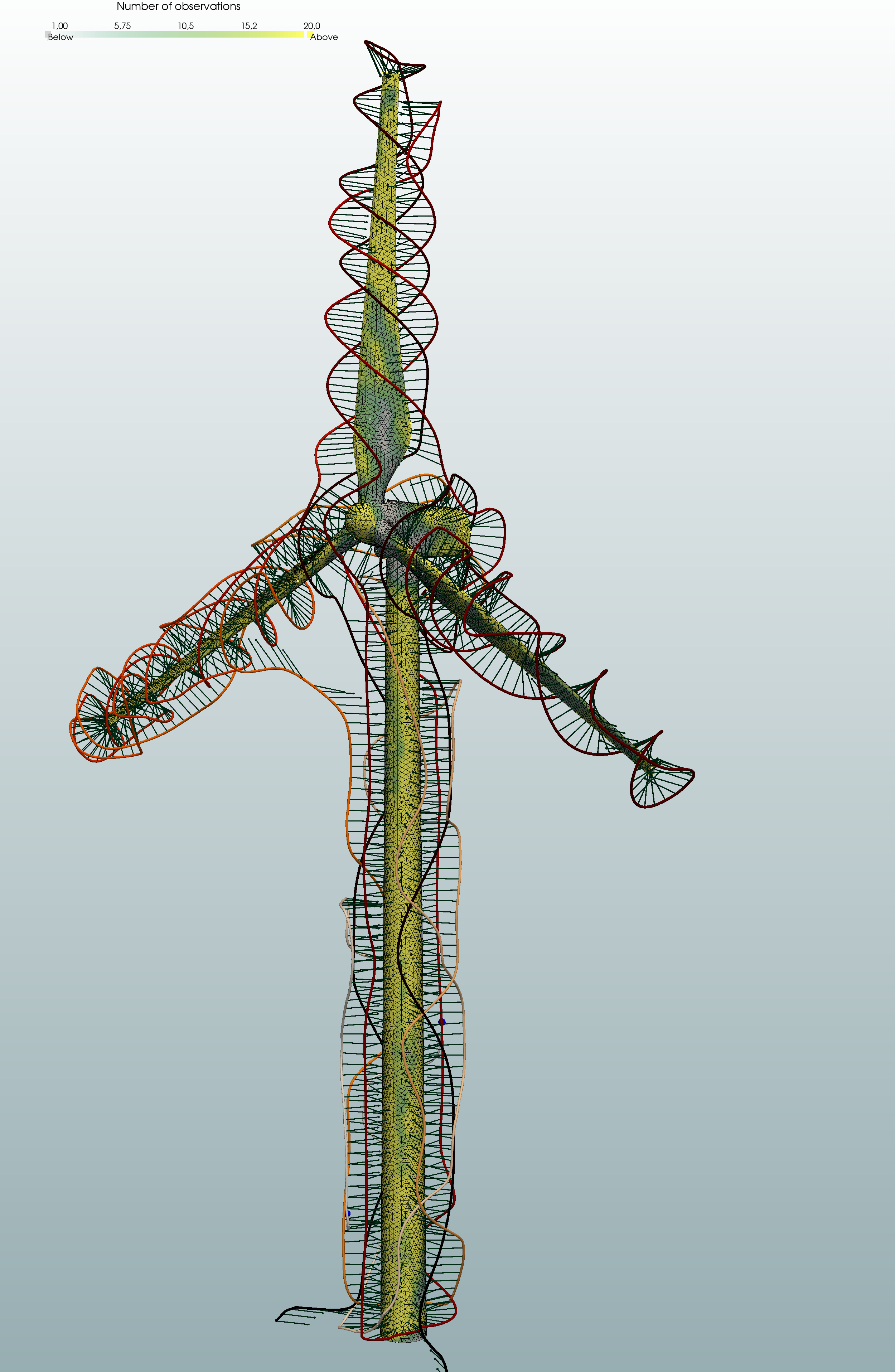}
	\caption{Relatively non-smooth convergence of the coverage measures $\eta_V$ and $\eta_A$ displayed in plot (A) is caused by alternating inspection motion between relatively isolated parts (column, hub, and blades) of the wind turbine structure. The plot of distances shows that, on average, the inspection distance is ensured and the minimal spacing constraint is complied. Trajectories of two UAVs generated by the HEDAC method for the inspection of the wind turbine in a total duration of 1 500 s are shown in (B). Camera orientations along trajectories are displayed as arrows of length equal to the inspection distance $d_m=6~\text{m}$. Green-to-yelow shading of the wind turbine surface indicates that almost all of it is observed between 5 to 20 times.}
	\label{fig:case_3_wind_turbine_inspection}
\end{figure}

A slightly stepped convergence of the spatial density coverage $\eta_{V}$ can be observed in Figure~\ref{fig:case_3_wind_turbine_inspection}(A) due to the UAVs' inspection focus switching between the wind turbine's column, hub, and blade. The Y-shaped configuration of the wind turbine is causing multi-pass inspections of the individual components, where each subsequent pass produces flatter change in $\eta_{V}$ and finally, it results in stepped coverage performance. Analogous to the first example, the distances of the trajectories from the structure are on average equal to $d_m=6~\text{m}$ and they have an increasing trend during the inspection. Since the trajectory planning is performed only for two UAVs and the domain is considerably large, the spacing between UAVs is easily kept at the safe distance $\epsilon$. 

Realized trajectories, camera orientations, and surface observations for wind turbine inspection are shown in Figure~\ref{fig:case_3_wind_turbine_inspection}(B). It is interesting to observe a synchronized spiral motion caused by the interaction between UAVs' coverage actions. All components of the turbine structure are observed from practically all directions and UAVs effectively collaborate in the inspection. The presented results imply that the proposed HEDAC trajectory planning algorithm is suitable for carrying out a multi-UAV inspection on real-world structures such as wind turbines.    

\subsection{Bridge test case}
Bridge inspection using UAVs is commonly conducted to visually detect damage as it lowers inspection costs and increases safety. Usually, UAVs are flown manually to visually inspect the structure using sensors and cameras. There has been a rising interest in research of UAV autonomous bridge inspection where previous studies have attempted to provide path planning and trajectory planning algorithms \cite{shi2021inspection}. The most noticeable one is the receding horizon next-best-view planner \cite{bircher2018receding} which is described in Section \ref{sec:two}. The authors employ the algorithm both for the exploration of a domain containing an unknown structure and for the visual inspection of a known structural model. For the inspection test case, the environment is represented as a volumetric occupancy grid map where segments are marked as inspected or uninspected based on the sensor's readings within 10 m, and paths are planned accordingly. To achieve greater coverage and improve computational performance, the planner produces paths by taking into account sensor readings at a distance of 2 m. A colliding box around the agent assures a safe distance between the agent and the structure. Visual inspection of a bridge structure using a single UAV lasted for $t_{tot}=105.0~\text{min}$ and achieved coverage of 99.1\%. Within this section, we reconstruct the same bridge test case, adapt it for HEDAC, and calculate the inspection trajectories with proper camera orientations. All parameters used for the bridge test case and their analogies from \cite{bircher2018receding} are provided in Table~\ref{tab:case_4_bridge_parameters}. Note that parameter values marked with $^b$ are used for comparison purposes as presented in subsection~\ref{subsec:comparison_receding_horizon}.

\begin{table}[!htb]
	\small
	\centering
	\caption{Geometry and numerical mesh info and parameters used for two variants of the bridge test case (values specific for inspection using five UAVs and single UAV are marked with $^a$ and $^b$ superscript, respectively) and corresponding or analogous parameters used in \cite{bircher2018receding} (marked with $^*$). Values marked with $^b$ are used in subsection~\ref{subsec:comparison_receding_horizon} where HEDAC is compared with Receding Horizon approach.}
	\begin{tabularx}{\linewidth}{Xrl}
		Parameter & Value & Unit\\
		\hline
		Structure length & 47.5 & m \\
		Structure width & 12.85 & m \\
		Structure height & 9.0 & m \\
		Domain length & 58.3, \scriptsize50$^*$ & m \\
		Domain width & 22.85, \scriptsize25$^*$ & m \\
		Domain height & 19.9, \scriptsize14$^*$ & m \\
		\hline
		Number of domain mesh nodes & 93 539 \\
		Number of domain mesh elements & 503 852 \\
		\scriptsize ~- Based on volumetric map resolution \cite{bircher2018receding} & \scriptsize 1 120 000$^*$\\
		Number of structure mesh surface nodes & 23 971 \\
		Number of structure mesh surface faces & 48 118 \\
		\scriptsize ~- Based on inspection mesh resolution \cite{bircher2018receding} & \scriptsize 140 000$^*$ \\
		\hline
		Inspection distance $d_m$ & 1.5 & m \\
		\scriptsize ~- $d_{planner\_max}$ \cite{bircher2018receding} & \scriptsize2$^*$& \scriptsize m \\
		Inspection distance broadness $d_\sigma$ & 0.3 & m \\
		\hline
		\hspace{0pt}FOV cone height $C_H$ & 2$^a$, 10$^b$ & m \\
		\hspace{0pt}FOV cone diameter $C_D$ & 3$^a$, 16$^b$ & m \\
		\hline
		Number of UAVs & 5$^a$, 1$^b$ \\
		UAV velocity $v$ & 0.5, \scriptsize 0.5$^*$ & m/s \\
		Safety distance $\epsilon$& 0.5 & m \\	
		\scriptsize ~- Colliding box \cite{bircher2018receding} & \scriptsize 0.5 $\times$ 0.5 $\times$ 0.3$^*$ & \scriptsize m \\
		Coverage action intensity $\Phi$ & 0.4$^a$, 1$^b$ \\
		Coverage action range $\sigma$ & 1 & m \\	
		\scriptsize ~- $d_{sensor\_max}$ \cite{bircher2018receding} & \scriptsize 10$^*$ & \scriptsize m \\ 
		\hline
		HEDAC conduction coefficient $k$ & 20 \\
		Inspection duration & 1 000$^a$, 3 000$^b$ & s \\
		Path planning time step $\Delta t$ & 0.5 & s \\
		\hline
	\end{tabularx}
	\label{tab:case_4_bridge_parameters}
\end{table}

\begin{figure}[!htb]
	\centering
	\includegraphics[width=\linewidth]{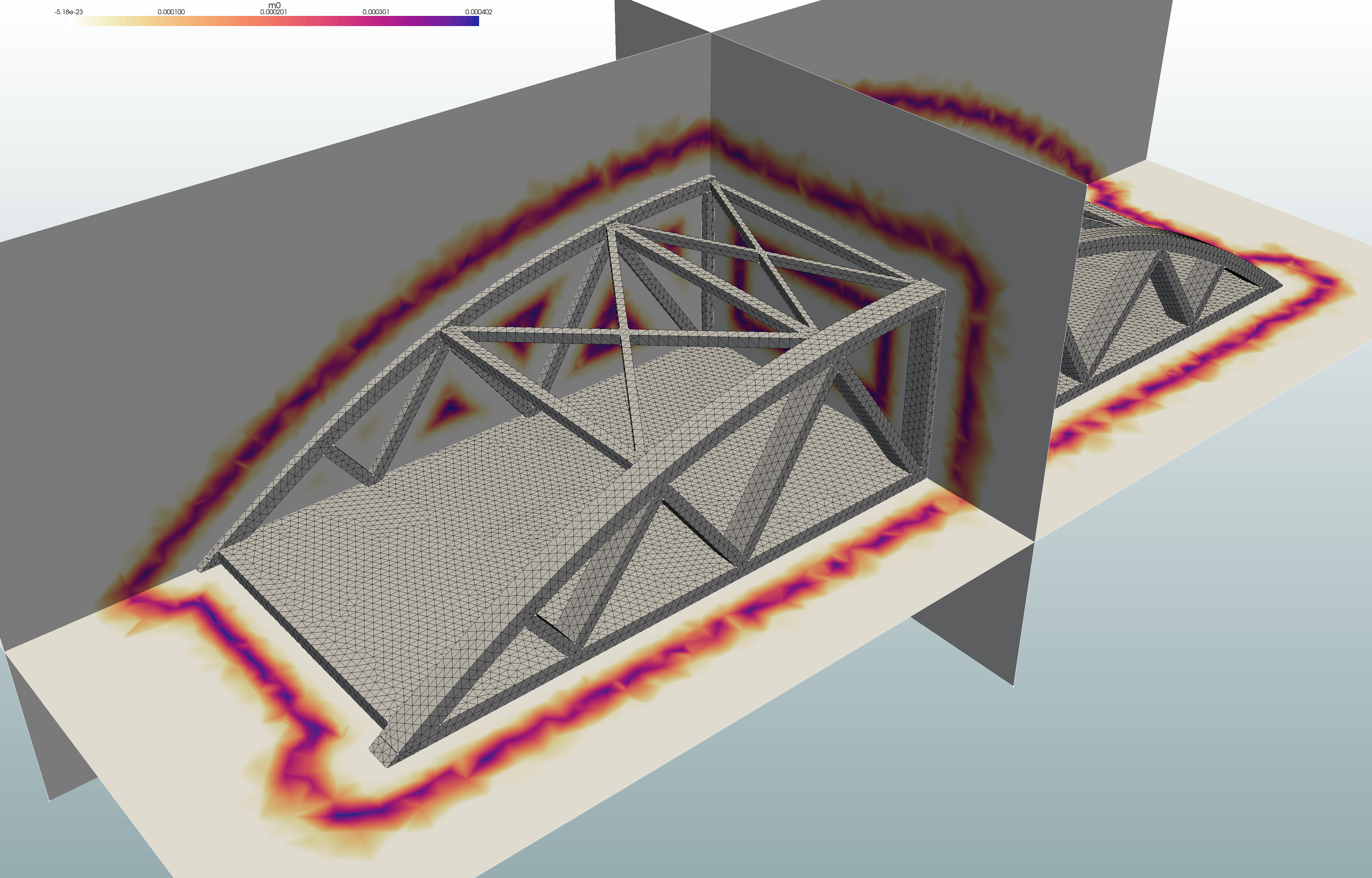}
	\caption{Spatial regions for bridge inspection trajectory planning are defined using a detailed target density field $\mu_0$. The field is constructed around the structure at the distance $d_m=1.5~\text{m}$ and with the inspection broadness $d_\sigma=0.3~\text{m}$.}
	\label{fig:case_4_bridge_mu0}
\end{figure}

The structure consists of around 50 m long and 13 m wide bridge deck as well as two 9 m high arches, increasing the inspection complexity when compared to the portal or wind turbine cases. The inspection domain is represented with the distance field inside the domain size volume. The distance field contains distances from the structure; therefore, leading agents near the structure to inspect. Based on the distance field $d_s$, we calculate the target density $\mu_0$ that covers the space around all structural components of the bridge at the distance $d_m$ with broadness $d_\sigma$ as shown in Figure~\ref{fig:case_4_bridge_mu0}. The implemented collision avoidance algorithm successfully prevents agent-to-agent as well as agent-to-boundary and agent-to-structure collisions. 
 
\begin{figure}[!htb]
	\centering	
	\includegraphics[width=\linewidth]{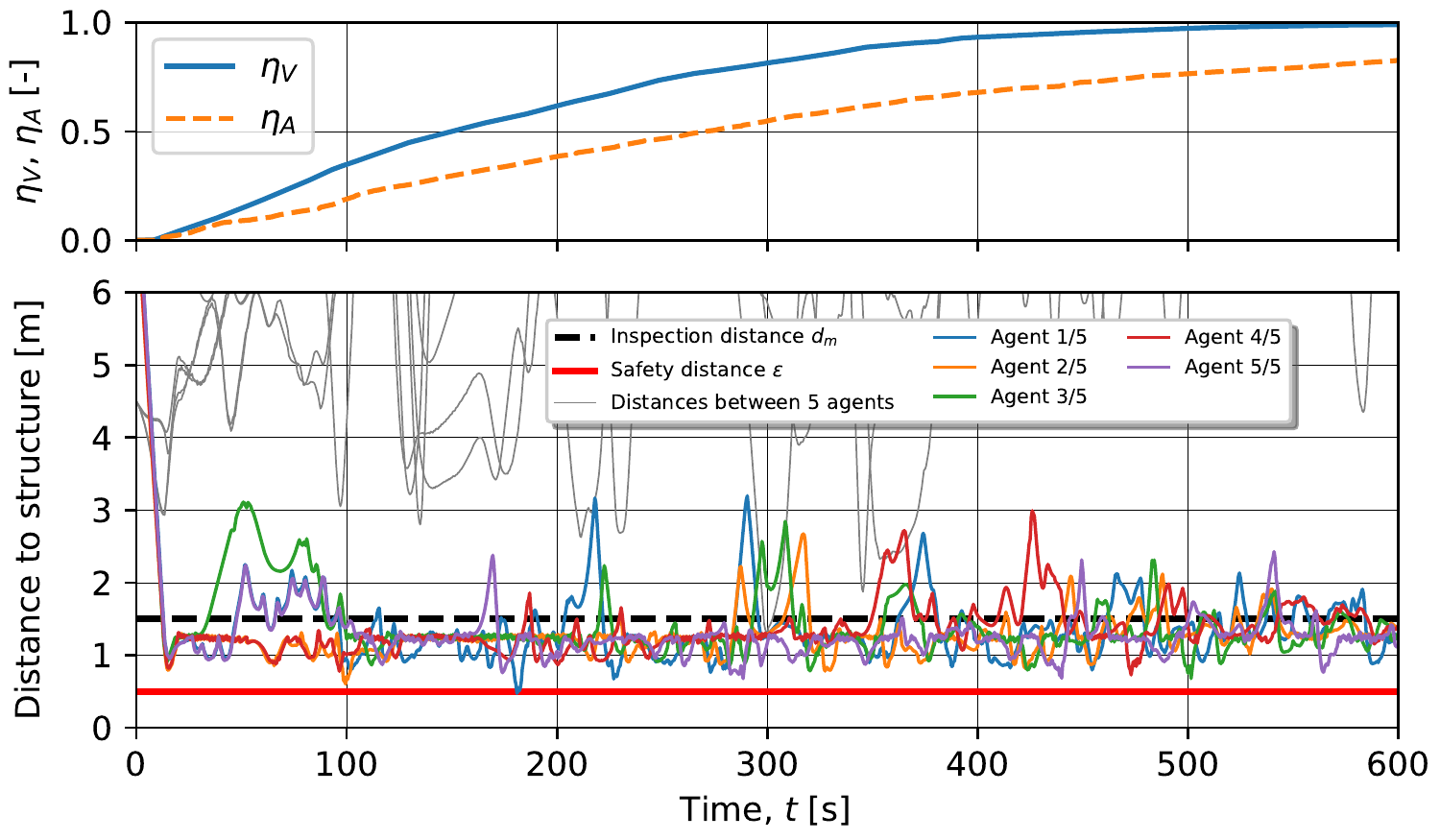}	
	\caption{Top subplot shows the convergence of spatial and surface coverage for the bridge inspection scenario using five UAVs. Calculated distances are shown in the bottom subplot. It is clearly visible that the inspection distance $d_m=1.5~\text{m}$ is hard to maintain during the operation due to the highly detailed and complex structure and flight domain. Although not as evident as in previous cases, the increasing trend of the UAV-to-structure distance is also present during the bridge inspection.}
	\label{fig:case_4_bridge_convergence}
\end{figure}

The convergence of the coverage measure $E$ and distances from agents to the structure are shown in Figure \ref{fig:case_4_bridge_convergence}. The inspection distance is difficult to maintain during the operation due to the complex structure and flight domain. However, the agents are restricted with the safety distance assuring collision avoidance. In a few brief moments, the safety distance is violated while inspecting the bridge with five agents due to the crowded flight domain, complex geometry, and large time step which prevents the collision avoidance from earlier activation. It can be avoided by reducing the time step or employing fewer agents for the inspection. 

\begin{figure}[!htb]
	\centering	
	\raggedright
	\includegraphics[width=\linewidth]{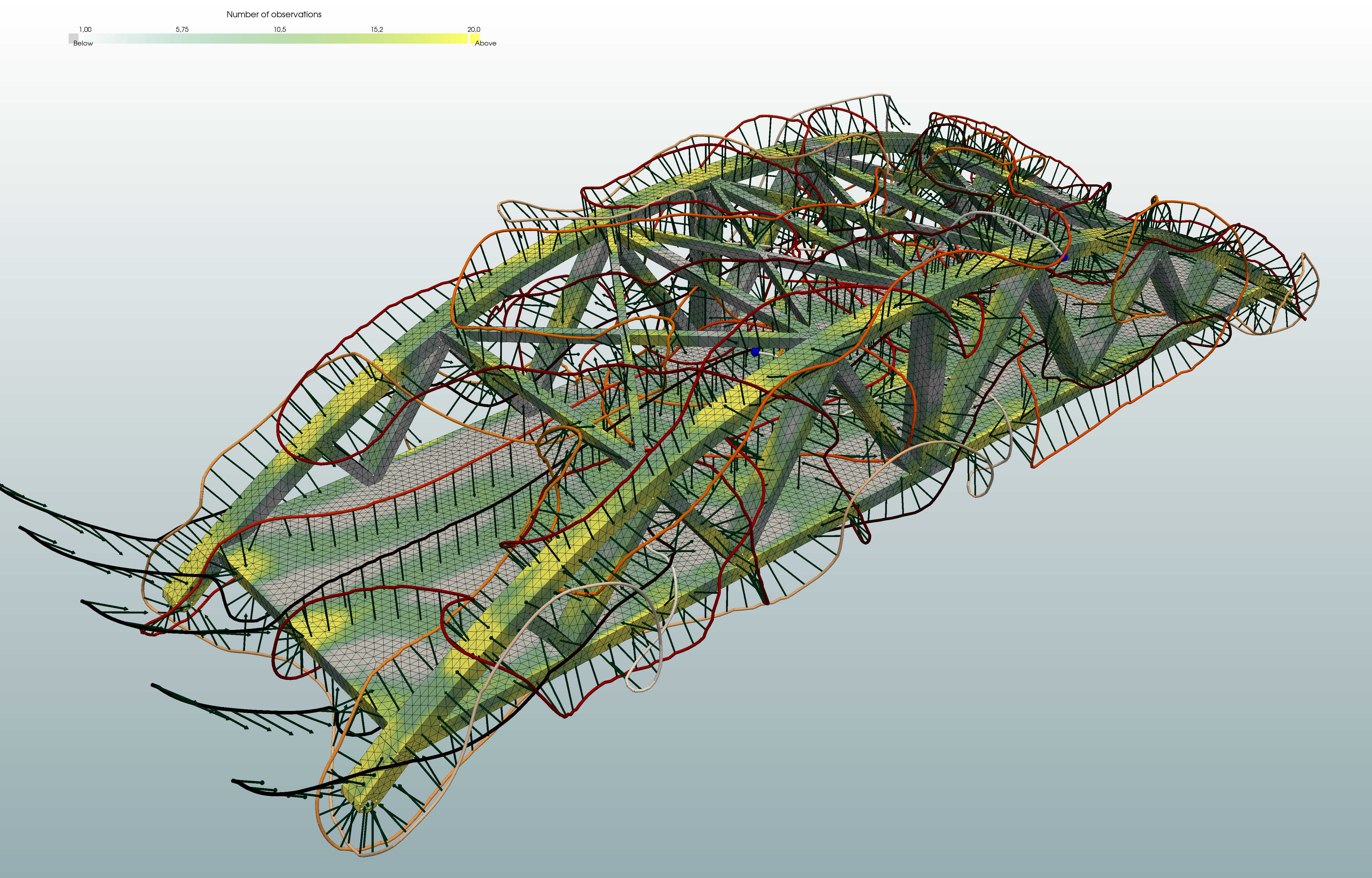}
	\raggedright 
	\caption{Obtained trajectories for five UAV bridge inspection planning. Trajectories successfully explore the domain and avoid frame trusses of the bridge structure. The directions of cameras are suitably adjusted and directed towards almost all surfaces of the bridge structure. Shading of the bridge surfaces represents the number of camera inspection shots achieved by movement of UAVs along generated trajectories. Since in the simulation each camera recording is triggered after each time step, the number of observations is accumulated. Due to the crossing of trajectories, the surfaces of some regions are observed more than 20 times. Remaining uninspected areas can be inspected by increasing the inspection time.}
	\label{fig:case_4_bridge_inspection}
\end{figure}

Produced trajectories are visualized in Figure \ref{fig:case_4_bridge_inspection}, showing appropriate camera orientations along the trajectories and achieved surface coverage.

%\newpage

\section{Performance and limitations analysis}
\label{sec:five}
The use of ergodic trajectory planning is well-suited for the UAV surface inspection tasks presented in this paper. Other approaches usually consider this as an optimization problem, which entails the potential of stalling in local optima and, in general, the limited effectiveness of the optimization itself.

Note that the parameter $k$ governs the balance between global (exploration) and local (exploitation) coverage. Due to the ergodic nature of the proposed method, it is self-balancing and both global and local coverage are eventually achieved. In theory, the perfect coverage is achieved in $\lim_{t \to \infty}$, while in practical applications, a near-optimal coverage can be achieved in a finite time. The duration of this finite time window depends on the choice of UAV motion and vision parameters ($\phi_\sigma$ or FOV) and on the scale and complexity of the inspected structure.

In the following subsections, we analyze the performance and shortcomings of the proposed method for UAV inspection trajectory planning. Validation of the method is carried out through comparison with a state-of-the-art approach on a standardized structure inspection scenario and is presented in subsection \ref*{subsec:comparison_receding_horizon}.

\subsection{Computational performance analysis}
%\added{The simulation of a UAVs path and camera orientation is simplified and does not represent faithfully real-world motion and time plan. Because of this limitation/simplification, this algorithm can be considered as path planner. For a more realistic model one should consider UAV characteristics and use a second order agent motion model which takes inertial effects into account.}
%\comment[Stefan]{Ovaj uvodni odlomak iznad mi nije dobar za ovo poglavlje. Možda se može iskoristiti negdje drugdje? A ovdje dodati jednu uvodnu rečenicu za computational benchmark.}
One of the main drawbacks of our proposed method is its computational inefficiency.
Handling a relatively large three-dimensional numerical mesh, when calculating coverage and potential, is very demanding. In this subsection, we analyze the computational costs of the method in more detail.

The initialization of the FEM in Table  \ref{tab:computational_efficiency} is computationally demanding because the sparse matrix is inverted, but it can still be done on a PC. This ensures a low-computational cost of FEM run-time calculation which reduces to only one matrix-vector multiplication at each time step.	
The optimization problem related to the collision avoidance can be computationally very demanding, but in this application it is formulated as an almost linear programming problem and does not require significant resources. 
The cost of the collision avoidance algorithm in practice can be seen in Table  \ref{tab:computational_efficiency}. The cost obviously depends on the particular test case, but it is almost insignificant compared to other parts of the algorithm.

%\begin{table}[!htb]
%	\small
%	\centering
%	\begin{tabularx}{\linewidth}{Xrrrrr}
%		Case &FEM init.&Coverage&FEM&Collision av.\\
%		\hline
%		Portal & 305s& 66.7\% & 33.2\%& 0.1\% \\
%		Turbine & 200s & 50.42\%  & 49.57\% & 0.01\% \\
%		Bridge & 130s& 54.6\% &44.4\%& 1\% \\
%		\hline
%	\end{tabularx}
%\caption{Execution time of FEM initialization and relative computational cost of coverage, FEM and collision avoidance calculation.}
%\label{tab:Execution time}
%\end{table}

\begin{table}[!htb]
	\footnotesize
	\centering
	\caption{Computational times for integral procedures of the proposed UAV trajectory planning algorithm. Procedures marked with $^*$ are being run in each time step and the displayed computational times are calculated as averages of all steps in each inspection simulation. 
	The computational benchmarks are executed using Python 3.8 on a Linux-based operating system, equipped with AMD Ryzen Threadripper 3970X 32-Core processor and 125,7 GiB of RAM.}
	\begin{tabularx}{\linewidth}{Xrrrrr}
		\hline
		~ & \multicolumn{3}{c}{UAV inspection scenario} \\
		Procedure 							& Portal 	& Wind turbine 	& Bridge \\
		\hline
		Initialization 						& 19.7 s 	& 114.2 s 		& 237.5 s  \\
		Coverage $\rho$ calculation$^*$ 	& 1.631 s 	& 7.415 s		& 5.132 s\\
		Potential $\psi$ calculation$^*$ 	& 0.832 s 	& 7.288 s		& 4.173 s\\
		Collision avoidance$^*$  			& 0.003 s	& 0.002 s		& 0.094 s\\
		Time step (all calculations)$^*$   	& 2.465 s	& 14.705 s		& 9.399 s \\	
		Entire trajectory planning 			& 3 717.2 s	& 35 406.2 s	& 19 035.5 s \\
		\hline
	\end{tabularx}
	\label{tab:computational_efficiency}
\end{table}

Although the proposed three-dimensional HEADC methodology is designed as a motion control algorithm, it is obvious that, due to excessive computational demands, it can not deliver results in real time. That is why we present this method as trajectory planning instead of motion control.

\subsection{Space coverage vs. surface coverage}

\begin{figure*}[!htb]
	\centering	
	\includegraphics[width=\linewidth]{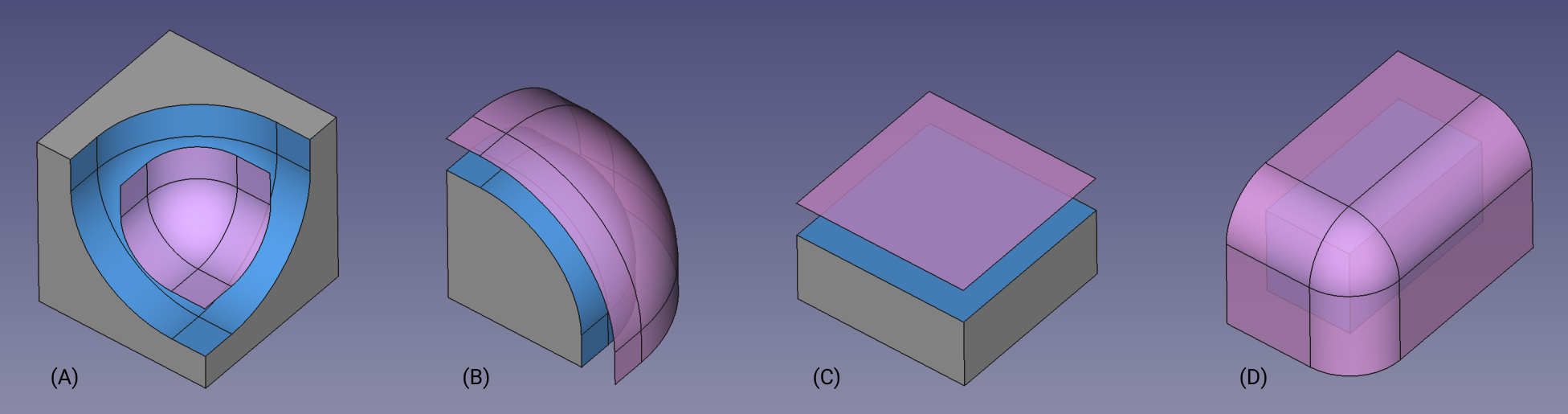}	
	\caption{The offset of a structure's surface depends on the surface curvature: (A) offset area of a concave surface decreases, (B) offset area of a convex surface increases, (C) offset area of a planar surface remains the same and (D) offset area for sharp edges and corners (of zero area) consists of cylindrical and spherical surfaces (infinite area increase). Referent and offset surfaces are marked with blue and transparent pink colors, respectively.}
	\label{fig:surfaces_offset}
\end{figure*}

The proposed trajectory planning method utilizes a space coverage technique in order to achieve coverage of two-dimensional surfaces. This approach is not justified in an obvious way, since space and surface coverage are not equivalent. However, as shown in the previous section, our approach can yield a suitable approximation, appropriate for real-world application. Nevertheless, the effect of this approximation needs to be addressed and assessed in pursuance of an overall evaluation of the proposed trajectory planning methodology.

The target density field $\mu_0$ is obtained by convolution of the Gaussian action $\phi_\sigma$ over the offsets of the structure's surfaces. The broadness of the target density field, regulated by the broadness $d_\sigma$, has an insignificant influence on the surface-to-space coverage mapping. However, the offset distance, regulated by the goal inspection distance $d_m$, plays a significant role. It is proportional to the surface-to-space mapping error. Using $d_\sigma=0$ (Gaussian function becomes Dirac function, providing infinitesimally thin goal density) and using offset $d_m=0$ would result in a spatial density $\mu_0$ that is equivalent to the structure's surface. Compared to the referent surface, the area of the offset surface can generally be the same, larger or smaller, for planar, convex, or concave surfaces, respectively (Figure~\ref{fig:surfaces_offset}). The influence of offsetting is not trivial \cite{zhuo2012curvature}, but it depends on surface curvature and offset distance.

\begin{figure}[!htb]
	\centering	
	\includegraphics[width=\linewidth]{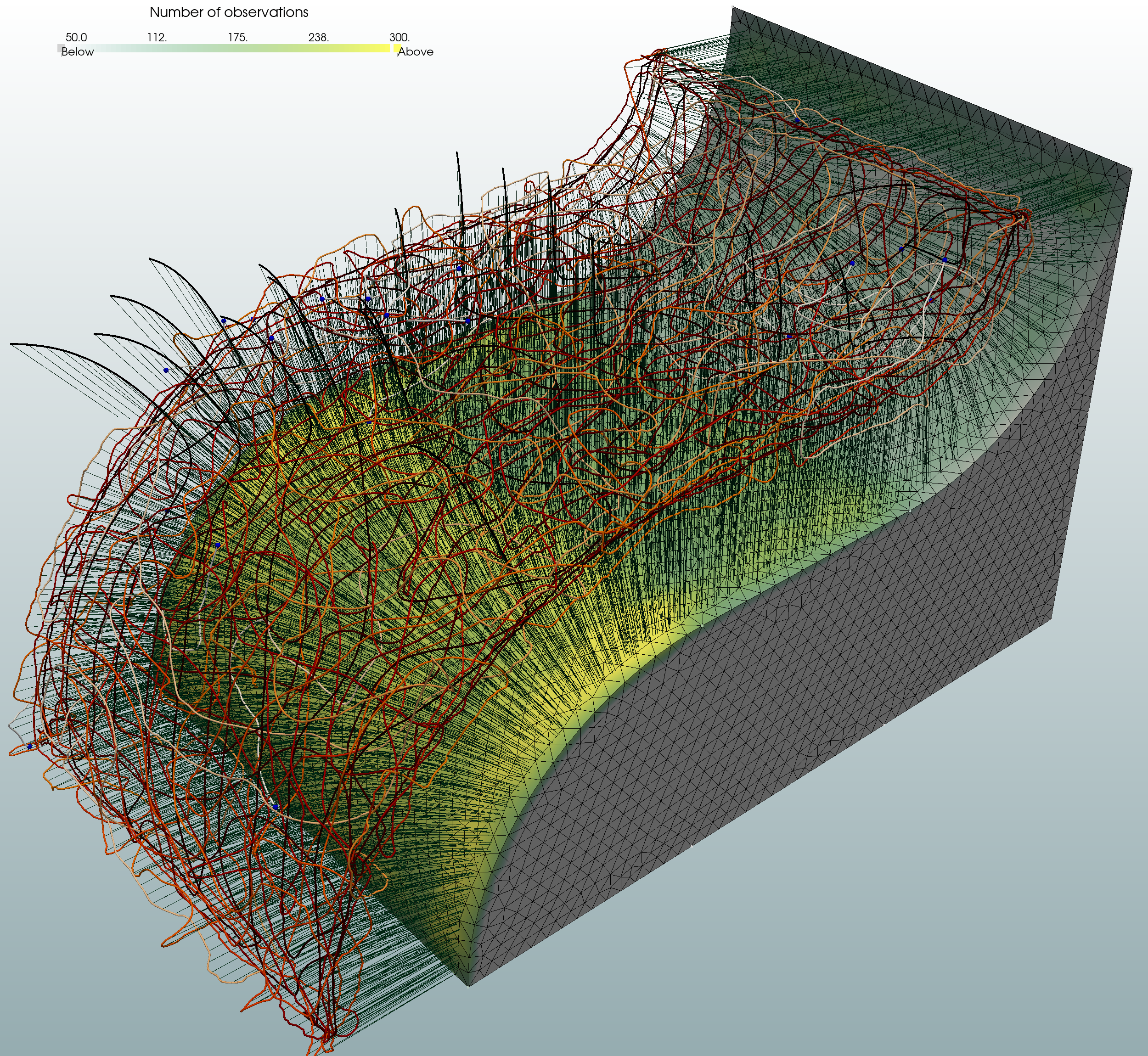}	
	\caption{Simulation of inspection of a surface with convex, concave and flat regions. All trajectories are located at the offset of the observed surface as defined by the target density $\mu_0$. Different numbers of recordings can be observed in each region of the plot, varying from 50 to 300. The non-linearity of delegating the surface coverage problem to space coverage is also indicated by the density of the camera direction vector (which coincides with the normal to the surface).}
	\label{fig:curved_surf_coverage}
\end{figure}

An inspection simulation on a relatively simple domain is prepared in order to demonstrate the effect of surface curvature on the inspection quality. The structure and the surrounding domain are designed to expose only one side of the S-shaped surface to the UAVs (Figure~\ref{fig:curved_surf_coverage}). The S-shaped surface has concave, convex, and flat regions. We simulate the inspection conducted by 20 UAVs for a duration of 1000 s. One can easily recognize a different number of achieved observations in each of the three regions, although the spatial target density $\mu_0$ is uniform (using the same distance $d_m$ and broadness $d_\sigma$) for all regions. These results confirm the considerations presented in the first part of this section.

Differences between spatial and surface coverage are strongly manifested in frequent and dense trajectories over certain regions.
Since a relatively small surface coverage is required for inspection (only a few shots of the same point are sufficient), achieving spatial coverage is an adequate solution for achieving surface coverage in UAV inspection tasks.

\subsection{Comparison with Receding Horizon}
\label{subsec:comparison_receding_horizon}

We utilize the receding horizon methodology for trajectory planning and the bridge test case, presented in \cite{bircher2018receding}, in order to evaluate the presented HEDAC trajectory generation approach. Since both compared approaches are designed differently, the parameters used in the trajectory planning procedures are different. In order to provide a scenario as equivalent as possible to the one in \cite{bircher2018receding}, we adjust the bridge case parameters (marked with $^b$ in Table~\ref{tab:case_4_bridge_parameters}). Additionally, non-identical but corresponding parameters used in the receding horizon approach are marked with $^*$.

\begin{figure}[!htb]
	\centering	
	\includegraphics[width=\linewidth]{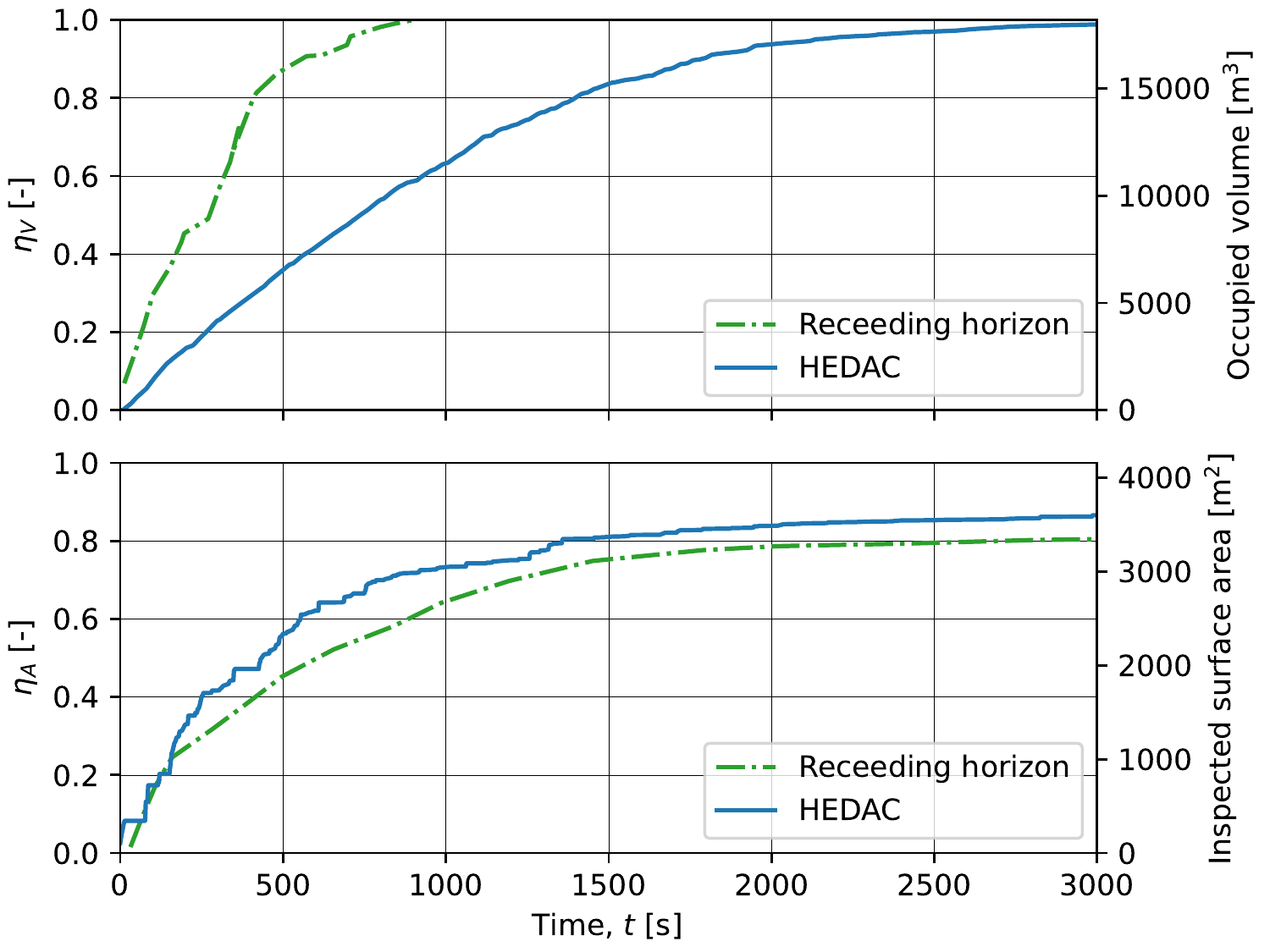}	
	\caption{Convergence of spatial coverage $\eta_V$ (or free volume \cite{bircher2018receding}, marked on the right-hand y-axis) for the bridge scenario using the proposed HEDAC method and receding horizon \cite{bircher2018receding} is shown on the upper plot. HEDAC occupies less volume by focusing on the space around the bridge surfaces allowing a more efficient surface inspection as shown in the lower plot. The lower plot shows surface coverage $\eta_A$ (or inspected surface area \cite{bircher2018receding}, marked on the right-hand y-axis).}
	\label{fig:case_4_bridge_comparison}
\end{figure}

The upper plot in Figure~\ref{fig:case_4_bridge_comparison} shows the comparison of spatial target density vs. free volume achieved with spatial ergodic coverage using the HEDAC method and volumetric exploration using the receding horizon \cite{bircher2018receding}, respectively. It is noticeable that the HEDAC approach occupies less space than the receding horizon. This can be easily explained by the fact that HEDAC uses only a part of the available space defined by $\mu_0>0$, while the receding horizon considers the entire volume available for flying. Focusing on the space around the surfaces being observed utilizes less volume and produces a more efficient inspection. This is visible in the lower plot of Figure~\ref{fig:case_4_bridge_comparison}, which shows the convergence of the area of inspected surfaces. Although it uses more conservative metrics (an effectively smaller FOV), the HEDAC approach outperforms the receding horizon for the bridge structure inspection presented in \cite{bircher2018receding}. 

\section{Conclusions}
\label{sec:six}
Due to the fast development of UAV technology and the potential cost benefits it could bring, autonomous tasks are prominent topics for contribution to the industries such as civil infrastructure inspection or search and rescue operations. To enable autonomous flight, safe paths need to be determined and provided to the UAVs assuring collision avoidance within the environment. In this paper, we have presented a new ergodic method for multiple UAV trajectory planning within a known three-dimensional environment as well as its application for visual infrastructure inspection. The method is an extension of the HEDAC algorithm, previously developed and tested for the exploration of two-dimensional domains. The presented algorithm calculates trajectories that cover three-dimensional space according to the given target density. The UAV motion is represented using a simple, first-order kinematic model allowing the calculation of trajectories feasible for a UAV. Collision avoidance is formulated as a non-linear optimization problem, successfully preventing collisions between UAVs as well as UAVs with domain boundaries by gradually redirecting a UAV when the distance threshold is reached. We have validated the implementation of collision avoidance in a crowded, unit cube domain with 100 UAVs uniformly covering the lower half of the cube. 

We have further adapted the algorithm for the visual inspection of three-dimensional structures using UAVs. In mathematical terms, the task is to explore two-dimensional manifolds embedded in a bounded three-dimensional volume using sensors mounted on a UAV, e.g., a camera. We have achieved this by constructing a three-dimensional field of interest (target density) as an offset from the structure's surface to be inspected. We have proposed a simple idea for handling the camera orientation. It relies on the gradient of the distance field $d_s$ containing distances to the nearest surface of the inspected structure. This solution produces camera orientations that always point to the nearest point on the inspected structure.

The proposed method has been tested on three inspection applications. The first test involves a synthetic portal scenario with 3 UAVs inspecting the structure. For a given target density field, continuous trajectories, safely distanced from the structure, are produced and the coverage convergence is achieved. The surfaces of the structure are covered relatively uniformly, leading to the conclusion that planned trajectories and camera orientations are suitable for infrastructure inspection with appropriate camera equipment. The second scenario of a wind turbine inspection provides a realistic test case demonstrating the coordination of 2 UAVs inspecting a relatively slender structure. As all components of the wind turbine are observed, we can implicatively conclude that the algorithm is suitable for conducting an autonomous UAV inspection of real-world structures. The third scenario of a bridge inspection has been successfully conducted with five UAVs inspecting the structure. 
The bridge scenario was also utilized in comparison with receding horizon inspection path planning. The results of the comparison indicate the competitiveness of the proposed approach, considering that the HEDAC method has achieved a better quality inspection.

Trajectories obtained with the proposed planner, though seemingly chaotic, are uniformly exploring the region around the 3D structure and allow the inspection of almost all surfaces of the structure in a limited time window. 
Although the differences between spatial and surface coverage are highlighted, the proposed approach is adequate for non-exhaustive, single-pass UAV inspection.
The convergence of the coverage measure is achieved for each scenario, with convergence depending on the complexity of the inspected structure, the number of UAVs inspecting the structure, and the available time window. The method is proven to be robust and stable, though it is not computationally efficient enough for real-time motion control applications. 

Some ideas for further research related to this topic, such as better domain mesh generation, adaptive mesh refinement, solving potential using General-Purpose computing on Graphics Processing Units (GPGPU), or domain partitioning, could potentially rise computational efficiency enough for real-time motion control implementation. Possible applications of the proposed methodology are numerous, not only for the tasks of UAV inspection of complex structures, but also for some contact actions such as UAV spraying or surface cleaning. A more accurate surface coverage solution could be solved by more advanced mapping of the surface to the target spatial density $\mu_0$ or by using a non-uniform motion model driven by the potential calculated directly on the surface of the structure using shell elements. The camera model is naive in the proposed inspection model (simplified FOV, zooming is not considered, and camera rotation is not limited), so for a real application, a suitable camera control should be designed in addition to UAV motion planning.

\section*{Acknowledgements}
This research is primarily supported by the Croatian Science Foundation under the project UIP-2020-02-5090. B.C.'s contribution is supported by Croatian Science Foundation under the project IP-2019-04-1239, L. M.’s contribution is supported by the European Union's Horizon 2020 Research and Innovation Programme under Grant Agreement No 861111, Drones4Safety.

\section*{Data availability}
All parameters for reproducing the study are presented in the manuscript. The data needed to reproduce the presented UAV inspection scenarios and video animations are available on the Open Science Framework repository: \url{{https://osf.io/bdrvn/}}. The Python code needed to reproduce this research is available upon request.

%\appendix
%\section{Sample Appendix Section}
%Lorem ipsum dolor sit amet, consectetur adipiscing elit, sed do eiusmod tempor section \ref{sec:sample1} incididunt ut labore et dolore magna aliqua. Ut enim ad minim veniam, quis nostrud exercitation ullamco laboris nisi ut aliquip ex ea commodo consequat. Duis aute irure dolor in reprehenderit in voluptate velit esse cillum dolore eu fugiat nulla pariatur. Excepteur sint occaecat cupidatat non proident, sunt in culpa qui officia deserunt mollit anim id est laborum.

 \bibliographystyle{elsarticle-num} 
 \bibliography{references}





% \end{thebibliography}
\end{document}